\documentclass{aims_noheader}
\usepackage{amsmath}
  \usepackage{paralist}
  \usepackage{graphics} 
\usepackage{subfigure}
  \usepackage{epsfig} 
\usepackage{graphicx}  \usepackage{epstopdf}
 \usepackage[colorlinks=true]{hyperref}
\hypersetup{urlcolor=blue, citecolor=red}

\def\currentvolume{}
 \def\currentissue{}
  \def\currentyear{}
   \def\currentmonth{}
    
\newtheorem{theorem}{Theorem}[section]

\theoremstyle{definition}

\newtheorem{example}[theorem]{Example}
\newtheorem{algorithm}[theorem]{Algorithm}

\def\barr{\begin{array}}
\def\earr{\end{array}}
\def\req#1{(\ref{#1})}
\def\norm#1{\hspace{0.2ex} \left\|#1\right\|}

\def\Tau{\tau}
\def\xsol{x^\dagger}

\def\ddd{\mathrm{d}}
\def\ttt{\mathrm{t}}
\def\tt{\tilde{\mathrm{t}}}
\def\rrr{\mathrm{r}}

\def\pp{{p}}
\def\pd{{{p}^*}}
\def\rr{{r}}

\def\ss{{s}}
\def\sd{{{s}^*}}

\def\Xd{{X^*}}

\def\xnd{x_n^\delta}
\def\xnpd{x_{n+1}^\delta}

\def\Djp#1#2{D_{\pp}(#1,#2)}
\def\Djpd#1#2{D_{\pd}(#1,#2)}
\def\Dpx0#1#2{D_\pp^{x_0}(#1,#2)}

\def\req#1{(\ref{#1})}
\def\N{{{\rm I}\!{\rm N}}}

\title[Parameters in a Newton- Landweber iteration in Banach space]{Enhanced choice of the parameters in an iteratively regularized Newton- Landweber iteration in Banach space}

\author[Barbara Kaltenbacher and Ivan Tomba]{}



 \email{barbara.kaltenbacher@aau.at}
 \email{ivan.tomba3@unibo.it}

\thanks{}

\begin{document}
\maketitle

\centerline{\scshape Barbara Kaltenbacher}
\medskip
{\footnotesize
 \centerline{Alpen-Adria-Universit\"at Klagenfurt}
   \centerline{Universit\"atstra\ss\ss e 65--67}
   \centerline{9020 Klagenfurt, Austria}
} 

\medskip

\centerline{\scshape Ivan Tomba}
\medskip
{\footnotesize
 \centerline{Universit\`{a} di Bologna}
   \centerline{Piazza Porta S. Donato, 5}
   \centerline{40127 - Bologna, Italy}
}

\bigskip

 \centerline{(Communicated by the associate editor name)}

\begin{abstract}
This paper is a close follow-up of \cite{KaltTomb13} and \cite{JinNLW}, where Newton-Landweber iterations have been shown to converge either (unconditionally) without rates or (under an additional regularity assumption) with rates. The choice of the parameters in the method were different in each of these two cases. We now found a unified and more general strategy for choosing these parameters that enables both convergence and convergence rates. Moreover, as opposed to the previous one, this choice yields strong convergence as the noise level tends to zero, also in the case of no additional regularity. Additionally, the resulting method appears to be more efficient than the one from \cite{KaltTomb13}, as our numerical tests show.
\end{abstract}

\section{Introduction}
Regularization of inverse problems in Banach spaces is a field of highly active research, cf., e.g.,
\cite{BurgerOsher:ConvergenceRatesVariationalRegularization,NHHKT10,scherzerpoeschl:convergencerates} for variational regularization and, e.g., \cite{BakuKok04,HeiKaz10,JinNLW,JinItTik,KalHof10,bkSchoepferSchuster09,SKHK12}, for iterative methods.
The main reason for this lies in the fact that extension of the scope from Hilbert to general Banach spaces better allows to formulate requirements on the searched for solution and to describe realistic noise models.

We will here especially concentrate on a combination of a Newton-type strategy with Landweber iterations to approximate the Newton step, which leads to a fully explicit iteration, cf. \cite{KaltTomb13,JinNLW}.

To formulate the method, consider a nonlinear ill-posed operator equation
\begin{equation}\label{Fxy}
F(x)=y
\end{equation}
where $F$ maps between Banach spaces $X$ and $Y$.
The given data $y^\delta$ are typically contaminated by noise, and we are going to assume that the noise level $\delta$ in
\begin{equation}\label{delta}
\|y-y^\delta\|\leq\delta
\end{equation}
is known.
In the following, $x_0$ is some initial guess and we will assume that a solution $\xsol$ to \req{Fxy} exists.
For some $\pp,\rr\in(1,\infty)$, we will make use of the duality mappings $J_\pp^X(x):=\partial \left\{ \frac{1}{\pp} \|x\|^\pp \right\}$ from $X$ to its dual $X^*$, and $J_\rr^Y(y):=\partial \left\{ \frac{1}{\rr} \|y\|^\rr \right\}$ from $Y$ to $Y^*$, respectively.
While under the assumptions we will make on $X$, the mapping $J_\pp^X$ will in fact be single valued, this will not necessarily be the case for $J_\rr^Y$ and we will denote by $j_\rr^Y$ a single valued selection from $J_\rr^Y$.
Therewith, we consider a combination of the iteratively regularized Gau\ss\ss-Newton method with an iteratively regularized Landweber method for approxi\-mating the Newton step, using some initial guess $x_0$ and starting from some $x_0^\delta$ (that need not necessarily coincide with $x_0$)
\begin{equation}\label{NIRLW}
\begin{array}{l}
\mbox{For }n=0,1,2\ldots\mbox{ do}\\
\hspace*{0.5cm}
u_{n,0}=0\\
\hspace*{0.5cm}
z_{n,0}=\xnd\\
\hspace*{0.5cm}
\mbox{For }k=0,1,2\ldots,k_n-1\mbox{ do}\\
\hspace*{1cm}
\begin{array}{l}
u_{n,k+1}= u_{n,k}-\alpha_{n,k}J_\pp^X(z_{n,k}-x_0) \\
\qquad\qquad -\omega_{n,k} F'(\xnd)^*j_\rr^Y (F'(\xnd)(z_{n,k}-\xnd ) + F(\xnd)-y^\delta)\\
z_{n,k+1}=x_0+{J_\pp^X}^{-1} \Bigl(J_\pp^X(\xnd -x_0) + u_{n,k+1}\Bigr)
\end{array}\\
\hspace*{0.5cm}
\xnpd=z_{n,k_n}
\,.
\end{array}
\end{equation}

It is clear that the choice of the parameters $\alpha_{n,k}$, $\omega_{n,k}$, $k_n$, and of the overall stopping index $n=n_*$ crucially influences the stability and efficiency of the method. While in \cite{KaltTomb13}, \cite{JinNLW} only either convergence or convergence rates have been established with disjoint parameter choices for each of these two cases,
the aim of this paper is to provide a unified parameter choice strategy for this method that allows to show both unconditional convergence and convergence rates under additional regularity assumptions on the solution.
Moreover, by using an appropriate choice of the stopping index for the inner iteration, differently from \cite{KaltTomb13}, \cite{JinNLW} we can show continuous dependence of the iterates $x_n$ on the data $y^\delta$ for each fixed outer iteration index $n$ and therewith are able to prove strong convergence.
Finally, the new parameter choice appears to enhance efficiency as compared to \cite{KaltTomb13}, as the numerical tests below show.

The remainder of this paper is organized as follows. In Section \ref{sec:prelim} we provide some preliminaries. The parameter choice as well as convergence results are derived and formulated in Section \ref{sec:conv_an}. Section \ref{sec:numres} shows some numerical tests for a coefficient identification problem in an elliptic PDE in one and two space dimensions.

\section{Preliminaries and assumptions}\label{sec:prelim}

Throughout this paper we will assume that $X$ is smooth, which means that the duality mapping is single-valued, and moreover, that $X$ is
$\ss$-convex for some $\ss\in[\pp,\infty)$, which implies
\begin{equation}\label{estcp}
\Djp{x}{y}\geq \frac{c_{p,s}}{(\|x\|+\|y\|)^{\ss-\pp}}\|x-y\|^{\ss}
\end{equation}
for some constant $c_{p,s}>0$, cf. Corollary 2.61 in \cite{SKHK12}.
Here, $\Djp{x}{y}$ denotes the Bregman distance
$$
\Djp{\tilde{x}}{x}=\frac{1}{\pp}\|\tilde{x}\|^\pp-\frac{1}{\pp}\|x\|^\pp
-\langle j_\pp^X(x),\tilde{x}-x\rangle_{X^*,X}
$$
(where $j_\pp^X(x)$ denotes a single valued selection of $J_\pp^X(x)$).
We will also make use of its shifted version
$$
\Dpx0{\tilde{x}}{x}:=\Djp{\tilde{x}-x_0}{x-x_0}\,.
$$
As a consequence of the above assumptions, $X$ is reflexive and we also have
\begin{equation}\label{estCq}
\hspace{-1cm}\Djpd{x^*}{y^*}\leq C_{p^*,s^*}\|x^*-y^*\|^{\sd}((p\Djpd{J_\pd^{X^*}(x^*)}{0})^{1-\frac{\sd}{\pd}}+\|x^*-y^*\|^{\pd-\sd})\,,
\end{equation}
for some $C_{p^*,s^*}$, where $\sd$ denotes the dual index $\sd=\frac{\ss}{\ss-1}$, cf. (4) in \cite{KaltTomb13}.
Under these assumptions,
the duality mapping is bijective and $J_\pp^{-1}=J_\pd^{X^*}$, the latter denoting the (by $\ss$-convexity also single-valued) duality mapping on the dual $X^*$ of $X$.
We will also make use of the identities
\begin{equation}\label{eq:three_point_identity}
\Djp{x}{y}=
\Djp{x}{z}+\Djp{z}{y}+\langle J_\pp^X(z)-J_\pp^X(y),x-z\rangle_{X^*,X}
\end{equation}
and
\begin{equation}\label{eq:connection_primal_dual_Bregman}
\Djp{y}{x}=\Djpd{J_\pp^X(x)}{J_\pp^X(y)}\,.
\end{equation}
For more details on the geometry of Banach spaces we refer, e.g., to \cite{schuster:iterative}, \cite{SKHK12} and the references therein.

The assumptions on the forward operator besides a condition on the domain
\begin{equation}\label{nonemptyint}
\mathcal{B}_{\rho}^D(\xsol)\subseteq \mathcal{D}(F)
\end{equation}
include a structural condition on its degree of nonlinearity. For simplicity of exposition we restrict ourselves to the tangential cone condition
\begin{equation}\label{tangcone}
\norm{F(\tilde{x})-F(x)-F'(x)(\tilde{x}-x)}
\leq\eta\norm{F(\tilde{x})-F(x)}\,, \
\tilde{x},x\in\mathcal{B}_{\rho}^D(\xsol)\,,
\end{equation}
and mention in passing that this could be extended to a more general condition on the degree of nonlinearity (cf. \cite{HeinHofm09}) as in \cite{IRLW}.
Here $F'$ is not necessarily the Fr\'{e}chet derivative of $F$ but just a linearization of $F$ s{at}isfying the Taylor remainder estimate \req{tangcone}.
Additionally, we assume that $F'$ and $F$ are uniformly bounded on $\mathcal{B}_{\rho}^D(\xsol)$.

Here
$$ \mathcal{B}_{\rho}^D(\xsol)=
\{x\in X\, | \,
\Dpx0{\xsol}{x}
\leq\rho^2\}
$$
with $\rho>0$ such that $x_0\in \mathcal{B}_{\rho}^D(\xsol)$.
By distinction between the cases $\|x-x_0\|<2\|\xsol-x_0\|$ and $\|x-x_0\|\geq 2\|\xsol-x_0\|$ and the second triangle inequality we obtain from \req{estcp} that
\begin{equation}\label{rhobar}
\mathcal{B}_{\rho}^D(\xsol)\subseteq \mathcal{B}_{\bar{\rho}}^{\|\cdot\|}(x_0)=
\{x\in X\, | \,
\|x-x_0\|
\leq\bar{\rho}\}
\end{equation}
with
$\bar{\rho}=\max\{ 2 \|\xsol-x_0\|\,, \ \left(\frac{2^\pp 3^{\ss-\pp}\rho^2}{c_{p,s}}\right)^{1/\pp}\}$

For obtaining convergence rates we impose a variational inequality
(or variational source condition)
\begin{eqnarray}\label{vinuNIRLW}
&&\exists\, \beta>0 : \ \forall x\in \mathcal{B}_{\rho}^D(\xsol) \quad \nonumber\\
&&|\langle J_\pp^X(\xsol-x_0),x-\xsol\rangle_{X^*\times X} |
\leq \beta \Dpx0{\xsol}{x}^{\frac{1-2\nu}{\ss}}
\|F(x)-F(\xsol)\|^{2\nu\frac{\lambda(\ss,\nu)}{\lambda(2,\nu)}} \,,
\end{eqnarray}
with with $s$ the parameter of smoothness of $X$, 
\begin{equation}\label{lambda}
 \lambda(\ss,\nu)=1-\frac{1-2\nu}{\ss}
\end{equation}
and $\nu\in(0,\frac12]$. 
Condition \eqref{vinuNIRLW} corresponds to a source condition $\xsol-x_0\in\mathcal{R}((F'(\xsol)^*F'(\xsol))^\nu)$ in the special case of Hilbert spaces (where $\ss=2$), cf., e.g., \cite{HeinHofm09}.
Note that \req{vinuNIRLW} is stronger for larger $\nu$ and \req{vinuNIRLW} always holds for $\nu=0$ with 
\begin{equation}\label{beta0}
\beta=\beta_0=\|\xsol-x_0\|^{\pp-1}\left(\frac{(\bar{\rho}+\|\xsol-x_0\|)^{\ss-\pp}}{c_{p,s}}\right)^{\frac{1}{s}}\,,
\end{equation}
due to \req{estcp}, \req{rhobar}.
The case $\nu=0$ can be identified with the situation that no additional regularity (i.e., \req{vinuNIRLW} with $\nu>0$) is known to hold.

\section{Error estimates and parameter choice}\label{sec:conv_an}
We here first of all follow the lines of Section 2 in \cite{KaltTomb13}. Some of the estimates are the same as in \cite{KaltTomb13} (and will be repeated here only for convenience of the reader). Some of them are different, though and therewith enable different parameter choice strategies.

For any $n\in\N$ we have
\begin{eqnarray}\label{estNIRLW1}
\lefteqn{\Dpx0{\xsol}{z_{n,k+1}}-\Dpx0{\xsol}{z_{n,k}}}\nonumber\\
&=&
\Dpx0{z_{n,k}}{z_{n,k+1}}
+\langle \underbrace{J_\pp^X(z_{n,k+1}-x_0)-J_\pp^X(z_{n,k}-x_0)}_{=
u_{n,k+1}-u_{n,k}},z_{n,k}-\xsol\rangle_{X^*\times X}\nonumber\\
&=& \underbrace{\Dpx0{z_{n,k}}{z_{n,k+1}}}_{(I)} \nonumber\\
&&-\underbrace{\omega_{n,k}\langle j_\rr^Y (A_n(z_{n,k}-\xnd)+F(\xnd ) - y^\delta), A_n (z_{n,k}-\xsol)\rangle_{Y^*\times Y}}_{(II)}\nonumber\\
&&-\underbrace{\alpha_{n,k} \langle J_\pp^X(\xsol-x_0),z_{n,k}-\xsol\rangle_{X^*\times X}}_{(III)}\nonumber\\
&&-\underbrace{\alpha_{n,k} \langle J_\pp^X(z_{n,k}-x_0)-J_\pp^X(\xsol-x_0), z_{n,k}-\xsol\rangle_{X^*\times X}}_{(IV)}\,,
\end{eqnarray}
where we abbreviate
\begin{equation}\label{Anbn}
A_n=F'(\xnd)\,,\quad b_n=y^\delta-F(\xnd)\,.
\end{equation}

Assuming that $z_{n,k}\in \mathcal{B}_{\rho}^D(\xsol)$,
we now estimate each of the terms on the right hand side separately.

By \req{estCq} and \req{eq:connection_primal_dual_Bregman} we have for the term (I)
\begin{equation}\label{estNIRLW_I}
\begin{aligned}
\lefteqn{\Dpx0{z_{n,k}}{z_{n,k+1}}}\\
\leq &C_{p^*,s^*}\|\underbrace{J_\pp^X(z_{n,k+1}-x_0)-J_\pp^X(z_{n,k}-x_0)}_{=u_{n,k+1}-u_{n,k}}\|^\sd\\
& \qquad \cdot\Bigl((p\rho^2)^{1-\frac{\sd}{\pd}}+
\|J_\pp^X(z_{n,k+1}-x_0)-J_\pp^X(z_{n,k}-x_0)\|^{\pd-\sd}\Bigr)\\
=& C_{p^*,s^*} (p\rho^2)^{1-\frac{\sd}{\pd}}\|\alpha_{n,k}J_\pp^X(z_{n,k}-x_0)+\omega_{n,k} A_n^*j_\rr^Y (
A_n(z_{n,k}-\xnd)+F(\xnd ) - y^\delta)\|^\sd\\
&+C_{p^*,s^*}\|\alpha_{n,k}J_\pp^X(z_{n,k}-x_0)+\omega_{n,k} A_n^*j_\rr^Y (
A_n(z_{n,k}-\xnd)+F(\xnd ) - y^\delta)\|^\pd\\
\leq& 2^{\sd-1}C_{p^*,s^*} (p\rho^2)^{1-\frac{\sd}{\pd}}
\alpha_{n,k}^\sd \|z_{n,k}-x_0\|^{(\pp-1)\sd}\\
&+2^{\sd-1}C_{p^*,s^*} (p\rho^2)^{1-\frac{\sd}{\pd}}
\omega_{n,k}^\sd \|A_n^*j_\rr^Y (A_n(z_{n,k}-\xnd)+F(\xnd ) - y^\delta)\|^\sd\Bigr)\\
&+2^{\pd-1}C_{p^*,s^*} \alpha_{n,k}^\pd
\|z_{n,k}-x_0\|^{(\pp-1)\pd}\\
&+2^{\pd-1}C_{p^*,s^*}
\omega_{n,k}^\pd \|A_n^*j_\rr^Y (A_n(z_{n,k}-\xnd)+F(\xnd ) - y^\delta)\|^\pd\Bigr)\\
\leq& C_{p^*,s^*} \Bigl(
(p\rho^2)^{1-\frac{\sd}{\pd}} \bar{\rho}^{(\pp-1)\sd} 2^{\sd-1} \alpha_{n,k}^\sd
+ \bar{\rho}^\pp 2^{\pd-1} \alpha_{n,k}^\pd\Bigr)
+ \varphi(\omega_{n,k} \tt_{n,k})
\end{aligned}
\end{equation}
where we have used the triangle inequality in $\Xd$ and $X$, the inequality
\begin{equation}\label{ablambda1}
(a+b)^\lambda\leq \left\{\begin{array}{ll}
2^{\lambda-1}(a^\lambda+b^\lambda)& \mbox{ if }\lambda\geq1\\
(a^\lambda+b^\lambda)& \mbox{ if }0\leq\lambda\leq1
\end{array}\right.\
\mbox{ for }a,b\geq0\,, \ \end{equation}
and \req{rhobar}, as well as the abbreviations
\begin{eqnarray}
\ddd_{n,k}&=&\Dpx0{\xsol}{z_{n,k}}^{1/2}\label{ddd}\\
\ttt_{n,k}&=&\|A_n(z_{n,k}-\xnd)+F(\xnd )-y^\delta\|\label{ttt}\\
\tt_{n,k}&=&\|A_n^*j_\rr^Y (A_n(z_{n,k}-\xnd)+F(\xnd )-y^\delta)\|\label{tt}\\
&\leq&\|A_n\|\ttt_{n,k}^{\rr-1}\,.\label{estttil}\\
\rrr_n&=&\|F(\xnd )-y^\delta\|\label{rrr}
\end{eqnarray}
Here
\begin{equation}\label{defvarphi}
\varphi(\lambda)= 2^{\sd-1}C_{p^*,s^*} (p\rho^2)^{1-\frac{\sd}{\pd}} \lambda^\sd + 2^{\pd-1}C_{p^*,s^*} \lambda^\pd\,,
\end{equation}
which by $\pd \geq \sd > 1$ defines a strictly monotonically increasing and convex function on $\mathbb{R}^+$.
Note that estimate \eqref{estNIRLW_I} is just the same as (15) in \cite{KaltTomb13}.

For the term (II) in \req{estNIRLW1} we get, using \req{tangcone}, \req{delta},
\begin{eqnarray}\label{estNIRLW_II}
\lefteqn{\omega_{n,k}\langle j_\rr^Y (A_n(z_{n,k}-\xnd)+F(\xnd ) - y^\delta), A_n (z_{n,k}-\xsol)\rangle_{Y^*\times Y}}\nonumber\\
&=&\omega_{n,k}\ttt_{n,k}^\rr \nonumber\\
&&+\omega_{n,k}\langle j_\rr^Y (A_n(z_{n,k}-\xnd)+F(\xnd ) - y^\delta),
\nonumber\\
&& \qquad \qquad  A_n (\xnd-\xsol)-F(\xnd ) +y^\delta\rangle_{Y^*\times Y}\nonumber\\
&\geq&\omega_{n,k}\ttt_{n,k}^\rr
-\omega_{n,k}\ttt_{n,k}^{\rr-1}(\eta\|F(\xnd ) - y^\delta\|+(1+\eta)\delta)
\end{eqnarray}
This is the same as (19) in \cite{KaltTomb13}.

To make use of the variational inequality \req{vinuNIRLW} for estimating (III), we first of all use \req{tangcone} to conclude
\begin{eqnarray*}
\lefteqn{\|F(z_{n,k})-F(\xsol)\|}\\
&=& \|(A_n(z_{n,k}-\xnd)+F(\xnd )-y^\delta)\\
&&\quad+ (F(z_{n,k})-F(\xnd)-A_n((z_{n,k}-\xnd)+(y^\delta-y)\|\\
&\leq&\ttt_{n,k}+\eta\|F(z_{n,k})-F(\xnd)\|+\delta\\
&\leq&\ttt_{n,k}+\eta(\|F(z_{n,k})-F(\xsol)\|+\|F(\xnd)-y^\delta\|)+(1+\eta)\delta\,,
\end{eqnarray*}
hence
\begin{equation}\label{FzFsol}
\|F(z_{n,k})-F(\xsol)\|\leq \frac{1}{1-\eta}\left(\ttt_{n,k}+\eta \rrr_n + (1+\eta)\delta\right)\,.
\end{equation}
This together with \req{vinuNIRLW} implies
\begin{eqnarray}\label{estNIRLW_IIIa}
|\lefteqn{\alpha_{n,k} \langle J_\pp^X(\xsol-x_0),z_{n,k}-\xsol\rangle_{X^*\times X}|}
\nonumber\\
&\leq&
\frac{\beta}{(1-\eta)^{2\nu}} \alpha_{n,k} \ddd_{n,k}^{1-2\nu}
\left(\ttt_{n,k}+\eta\rrr_n+(1+\eta)\delta\right)^{2\nu}
\nonumber\\
&\leq& C(\lambda(\ss,\nu))\frac{\beta}{(1-\eta)^{2\nu}} \alpha_{n,k}
\left\{ \ddd_{n,k}^2 + \left(\ttt_{n,k}+\eta\rrr_n+(1+\eta)\delta\right)^{\frac{4\nu}{1+2\nu}}\right\}
\end{eqnarray}
where we have used the elementary estimate
\begin{equation}\label{ablambda2}
a^{1-\lambda}b^\lambda\leq C(\lambda)(a+b)
\ \mbox{ for }a,b\geq0\,, \ \lambda\in(0,1)
\end{equation}
with $C(\lambda)=\lambda^\lambda(1-\lambda)^{1-\lambda}$ and $\lambda(\ss,\nu)$ as in \eqref{lambda}.
Note that \eqref{estNIRLW_IIIa} differs from the corresponding estimate (24) in \cite{KaltTomb13}.

Finally, for the term (IV) we have that
\begin{eqnarray}\label{estNIRLW IV}
\lefteqn{\alpha_{n,k}\langle J_\pp^X(\xsol-x_0)-J_\pp^X(z_{n,k}-x_0),\xsol-z_{n,k}\rangle_{X^*\times X}}\nonumber\\
&=&\alpha_{n,k}(\Dpx0{\xsol}{z_{n,k}}+\Dpx0{z_{n,k}}{\xsol})
\ \geq \ \alpha_{n,k} \ddd_{n,k}^2\,,
\end{eqnarray}
which is just (26) in \cite{KaltTomb13}.

Altogether we arrive at the estimate
\begin{eqnarray} \label{estab}
\ddd_{n,k+1}^2&\leq&
\Bigl(1-(1-c_0)\alpha_{n,k}\Bigr)\ddd_{n,k}^2
+c_1\alpha_{n,k}^\sd
+c_2\alpha_{n,k}^\pd
\nonumber\\
&&
+c_3 \alpha_{n,k}
\left(\ttt_{n,k}+\eta\rrr_n+(1+\eta)\delta\right)^{\frac{4\nu}{1+2\nu}}
\nonumber\\
&& -\omega_{n,k}\ttt_{n,k}^\rr + \omega_{n,k}\ttt_{n,k}^{\rr-1}(\eta \rrr_n + (1+\eta)\delta)
+ \varphi(\omega_{n,k} \tt_{n,k})\,,
\end{eqnarray}
where
\begin{eqnarray}
c_0&=&
\frac{\beta}{(1-\eta)^{2\nu}} C(\lambda(\ss,\nu))
\label{c0}\\
c_1&=&C_{p^*,s^*} (p\rho^2)^{1-\frac{\sd}{\pd}} \bar{\rho}^{(\pp-1)\sd} 2^{\sd-1}
\label{c1}\\
c_2&=&C_{p^*,s^*} \bar{\rho}^{\pp} 2^{\pd-1}
\label{c2}\\
c_3&=&\frac{\beta}{(1-\eta)^{2\nu}} C(\lambda(\ss,\nu))
\label{c3}\\
\theta&=&\frac{4\nu}{\rr(1+2\nu)-4\nu}\,, \label{theta}
\end{eqnarray}
(small $c$ denoting constants that can be made small by assuming $x_0$ to be sufficiently close to $\xsol$ and therewith $\beta$, $\eta$, $\|x_0-\xsol\|$ small), and $\lambda(\ss,\nu)$
as in \eqref{lambda}.

Multiplying \req{estab} with $\alpha_{n,k+1}^{-\theta}$ and abbreviating
$$
\gamma_{n,k}:=\ddd_{n,k}^2
\alpha_{n,k}^{-\theta}\,,
$$
we get
\begin{equation}\label{estgamma1}
\begin{aligned}
\lefteqn{\gamma_{n,k+1}-\gamma_{n,k}}\\
\leq \left(\frac{\alpha_{n,k}}{\alpha_{n,k+1}}\right)^\theta&
\Bigl\{\Bigl(  1-(1-c_0)\alpha_{n,k}-\left(\frac{\alpha_{n,k+1}}{\alpha_{n,k}}\right)^\theta
\Bigr)\gamma_{n,k}\\
&+(
c_1\alpha_{n,k}^{\sd-\theta}+c_2\alpha_{n,k}^{\pd-\theta}
+c_3 \alpha_{n,k}^{1-\theta} \left(\ttt_{n,k}+\eta\rrr_n+(1+\eta)\delta\right)^{\frac{4\nu}{1+2\nu}}
\\
& - \alpha_{n,k}^{-\theta}\left(\omega_{n,k}\ttt_{n,k}^\rr - \omega_{n,k}\ttt_{n,k}^{\rr-1}(\eta \rrr_n + (1+\eta)\delta) - \varphi(\omega_{n,k} \tt_{n,k})\right)\Bigr\}\,.
\end{aligned}
\end{equation}

To obtain monotone decay of the sequence $\gamma_{n,k}$ with increasing $k$ we choose
\begin{itemize}
\item $\omega_{n,k}\geq0$ such that
\begin{equation}\label{omegank1}
\underline{\omega}\leq\omega_{n,k}\leq\overline{\omega} \mbox{ and }
\frac{\varphi\Bigl(\omega_{n,k}\tt_{n,k}\Bigr)}{\omega_{n,k}\ttt_{n,k}^\rr}
\leq \overline{c}_\omega
\end{equation}
for some $0<\underline{\omega}<\overline{\omega}$, $\overline{c}_\omega>0$.
We will do so by setting
\begin{equation}\label{omega_explicit}
\omega_{n,k} =\vartheta\min\{
\ttt_{n,k}^{\frac{\rr}{\sd-1}}\tt_{n,k}^{-s}\,, \
\ttt_{n,k}^{\frac{\rr}{\pd-1}}\tt_{n,k}^{-p}\,, \ \overline{\omega}\}
\end{equation}
with $\vartheta$ sufficiently small, cf. \cite{KaltTomb13}, and assuming that
\begin{equation}\label{pnu}
\rr\geq \ss\geq \pp\,,
\end{equation}

\item $\alpha_{n,k}\geq0$ such that
\begin{equation}\label{alphank0}
\alpha_{n,k}\geq \check{\alpha}_{n,k}:=
\tilde{\tau}\left(\ttt_{n,k}+\eta\rrr_n+(1+\eta)\delta\right)^{\frac{\rr}{1+\theta}}
\end{equation}
and
\begin{equation}\label{alphank1}
c_0
+\frac{c_1}{\overline{\gamma}_{0,0}}\alpha_{n,k}^{\sd-\theta-1}+\frac{c_2}{\overline{\gamma}_{0,0}}\alpha_{n,k}^{\pd-\theta-1}+\frac{c_3}{\tilde{\tau}^\theta\overline{\gamma}_{0,0}} +\frac{1}{\tilde{\tau}^{1+\theta}\overline{\gamma}_{0,0}}\leq q<1\,.
\end{equation}
The latter can be achieved by
\begin{eqnarray}
\alpha_{n,k}\leq 1\mbox{ and }\label{alphank1_suff}\\
\sd\geq\theta+1 \,, \quad \pd\geq\theta+1\,, \label{rp1}\\
c_0,
c_1,c_2,c_3,\tilde{\tau}\mbox{ sufficiently small.} \nonumber
\end{eqnarray}
In case $\nu>0$ we additionally require
\begin{equation}\label{alphank2}
\alpha_{n,k+1}\geq\hat{\alpha}_{n,k+1}:=\alpha_{n,k} \Bigl( 1-(1-q)\alpha_{n,k}\Bigr)^{1/\theta}
\end{equation}
with an upper bound $\overline{\gamma}_{0,0}$ for $\gamma_{0,0}$. Note that this just means $\alpha_{n,k+1}\geq0$ in case  $\nu=0$, i.e., $\theta=0$, thus an empty condition in this case.

To meet conditions \req{alphank0}, \req{alphank2} with a minimal $\alpha_{n,k+1}$ we set
\begin{eqnarray}\label{alphank02}
&&\alpha_{n,k+1}=\max\{\check{\alpha}_{n,k+1}\,, \,
\hat{\alpha}_{n,k+1}\} \mbox{ for } k\geq 0 \\
&&\alpha_{n,0}=\left\{\begin{array}{ll}\alpha_{n-1,k_{n-1}}&\mbox{if }n\geq1\\\alpha_{0,0}&\mbox{if }n=0\end{array}\right.\,.\nonumber
\end{eqnarray}
\end{itemize}
It remains to choose
\begin{itemize}
\item
the inner stopping index $k_n$ and
\item the outer stopping index $n_*$,
\end{itemize}
see below.\\
Indeed with these choices of $\omega_{n,k}$ and $\alpha_{n,k+1}$ we can inductively conclude from \req{estgamma1} that
\begin{eqnarray}\label{mon}
\gamma_{n,k+1}-\gamma_{n,k}&\leq& \left(\frac{\alpha_{n,k}}{\alpha_{n,k+1}}\right)^\theta
\Bigl\{\Bigl(  1-(1-q)\alpha_{n,k}-\left(\frac{\alpha_{n,k+1}}{\alpha_{n,k}}\right)^\theta
\Bigr)\overline{\gamma}_{0,0}\Bigr\}
\nonumber\\&&
- \alpha_{n,k+1}^{-\theta}(1-\overline{c}_\omega)\omega_{n,k}\ttt_{n,k}^\rr\,,
\nonumber\\
&\leq& - \alpha_{n,k+1}^{-\theta}(1-\overline{c}_\omega)\omega_{n,k}\ttt_{n,k}^\rr\leq 0\,.
\end{eqnarray}
This monotonicity result holds for all $n\in\N$ and for all $k\in\N$.

By \req{mon} and $\alpha_{n,k}\leq 1$ (cf. \req{alphank1_suff}) it can be shown inductively that all iterates remain in $\mathcal{B}_{\rho}^D(\xsol)$ provided
\begin{equation}\label{gammarho}
\overline{\gamma}_{0,0}\leq \rho^2\,.
\end{equation}

Moreover, \req{mon} implies that
\begin{equation}\label{sumsumtnk}
\sum_{n=0}^\infty \sum_{k=0}^\infty\alpha_{n,k+1}^{-\theta}\omega_{n,k}\ttt_{n,k}^\rr \leq
\frac{\overline{\gamma}_{0,0}}{1-\overline{c}_\omega}<\infty,
\end{equation}
hence by $\alpha_{n,k+1}\leq 1$, $\omega_{n,k}\geq \underline{\omega}$
\begin{equation}\label{tnkto0k}
\ttt_{n,k}\to 0 \mbox{ as }k\to\infty \mbox{ for all }n\in\N
\end{equation}
and
\begin{equation}\label{tnkto0n}
\sup_{k\in\N_0}\ttt_{n,k}\to 0 \mbox{ as }n\to\infty\,.
\end{equation}
Especially, since $\ttt_{n,0}=\rrr_n$, \eqref{sumsumtnk} implies
\begin{equation}\label{rnto0}
\sum_{n=0}^\infty \rrr_n^\rr < \infty,
\end{equation}
hence $\rrr_n\to 0$  as $n\to\infty$.

To quantify the behavior as $k\to\infty$ of the sequence $\alpha_{n,k}$ according to \req{alphank0}, \req{alphank2}, \req{alphank02} for fixed $n$ we distinguish between two cases.\\
\textbf{Case (i)}: there exists a $\underline{k}$ such that for all $k\geq \underline{k}$ we have $\alpha_{n,k}=\hat{\alpha}_{n,k}$.
Considering an arbitrary accumulation point $\bar{\alpha}_n$ of $\alpha_{n,k}$ (which exists since $0\leq\alpha_{n,k}\leq 1$) we therefore have $\bar{\alpha}_n=\bar{\alpha}_n\Bigl(1-(1-q)\bar{\alpha}_n\Bigr)^\frac{1}{\theta}$, hence $\bar{\alpha}_n=0$.\\
\textbf{Case (ii)}: consider the situation that (i) does not hold, i.e., there exists a subsequence $k_j$ such that for all $j\in\N$ we have $\alpha_{n,k_j}=\check{\alpha}_{n,k_j}$. Then by \req{alphank0}, \req{alphank2}, and \req{tnkto0k} we have $\alpha_{n,k_j}\to
\tilde{\tau}\left(\eta\rrr_n+(1+\eta)\delta\right)^{\frac{r}{1+\theta}}$.\\
Altogether we have shown that
\begin{equation}\label{limsupalphak}
\limsup_{k\to\infty}\alpha_{n,k}\leq
\tilde{\tau}\left(\eta\rrr_n+(1+\eta)\delta\right)^{\frac{r}{1+\theta}}\mbox{ for all }n\in\N\,.
\end{equation}
Since $\eta$ and $\delta$ can be assumed to be sufficiently small, this especially implies the bound $\alpha_{n,k}\leq1$ in \req{alphank1_suff}.

We consider $z_{n_*,k_{n_*}^*}$ as our regularized solution, where $n_*$, $k_{n_*}^*$ (and also $k_n$ for all $n\leq n_*-1$; note that $k_{n_*}^*$ is to be distinguished from $k_{n_*}$ - actually the latter is not defined, since we only define $k_n$ for $n\leq n_*-1$!) are still to be chosen appropriately, according to the requirements from the proofs of
\begin{itemize}
\item convergence rates in case $\nu,\theta>0$,
\item convergence for exact data $\delta=0$,
\item convergence for noisy data as $\delta\to0$.
\end{itemize}

\subsection{Convergence rates in case $\nu,\theta>0$}\label{subsec:rates}
From \req{mon} we get
\begin{equation}\label{ratealpha}
\ddd_{n,k}^2\leq \overline{\gamma}_{0,0}\alpha_{n,k}^\theta \mbox{ for all }n,k\in\N\,,
\end{equation}
hence in order to get the desired rate
$$ \ddd_{n_*,k_{n_*}^*}^2=O(\delta^{\frac{r\theta}{1+\theta}})$$
in view of \req{limsupalphak} (which is a sharp bound in case (ii) above) we need to have a bound
\begin{equation}\label{discrprinc_upper}
\rrr_{n*}\leq \Tau\delta
\end{equation}
for some constant $\Tau>0$, and we should choose $k_{n_*}^*$ large enough so that
\begin{equation}\label{alphanstar}
\alpha_{n_*,k_{n_*}^*}\leq C_\alpha(\rrr_{n_*}+\delta)^{\frac{r}{1+\theta}}
\end{equation}
which is possible with a finite $k_{n_*}^*$ by \req{limsupalphak} for $C_\alpha>\left(\tilde{\tau}(1+\eta)\right)^{\frac{r}{1+\theta}}$.
Note that this holds without any requirements on $k_n$ for $n<n_*$.

\subsection{Convergence as $n\to\infty$ for exact data $\delta=0$}\label{subsec:conv_ex}


To show that $(x_n)_{n\in\N}$ is a Cauchy sequence (following the seminal paper \cite{HaNeSc95}), for arbitrary $m<j$, we choose the index $l\in\{m,\ldots,j\}$ such that $\rrr_l$ is minimal and use the identity
\begin{eqnarray}\label{eq:id}
\Dpx0{x_l}{x_m} &=& \Dpx0{\xsol}{x_m}-\Dpx0{\xsol}{x_l}\nonumber\\
&&+\langle J_\pp^X(x_l-x_0)-J_\pp^X(x_m-x_0),x_l-\xsol\rangle_{X^*\times X}
\end{eqnarray}
and the fact that the monotone decrease and boundedness from below
of the sequence $\Dpx0{\xsol}{x_m}$ implies its convergence, hence
it suffices to prove that the last term in \req{eq:id} tends to zero as $m<l\to\infty$. (Analogously it can be shown that $\Dpx0{x_l}{x_j}$ tends to zero as $l<j\to\infty$).
This term can be rewritten as
\begin{eqnarray*}
\lefteqn{\langle J_\pp^X(x_l-x_0)-J_\pp^X(x_m-x_0),x_l-\xsol\rangle_{X^*\times X}}\\
&=& \sum_{n=m}^{l-1} \sum_{k=0}^{k_n-1}
\langle u_{n,k+1}-u_{n,k},x_l-\xsol\rangle_{X^*\times X}\,,
\end{eqnarray*}
where
\begin{eqnarray*}
\lefteqn{\vert \langle u_{n,k+1}-u_{n,k} ,x_l-\xsol\rangle_{X^*\times X}\vert}\\
&=& \vert \alpha_{n,k} \langle J_\pp^X(z_{n,k}-x_0),
x_l-\xsol\rangle_{X^*\times X}\\
&&+\omega_{n,k} \langle j_\rr^Y (A_n(z_{n,k}-\xnd ) - b_n),
A_n (x_l-\xsol)\rangle_{X^*\times X} \vert \\
&\leq& 2\bar{\rho}^\pp \alpha_{n,k}
+\omega_{n,k} \ttt_{n,k}^{\rr-1} \|A_n (x_l-\xsol)\|  \\
&\leq& 2\bar{\rho}^\pp\tilde{\tau}(\ttt_{n,k}+\eta\rrr_n)^\rr+\omega_{n,k} \ttt_{n,k}^{\rr-1} (1+\eta)(2\rrr_n+\rrr_l)\\
&\leq& 2\bar{\rho}^\pp\tilde{\tau}(\ttt_{n,k}+\eta\rrr_n)^\rr+3(1+\eta)\omega_{n,k} \ttt_{n,k}^{\rr-1} \rrr_n\\
\end{eqnarray*}
by our choice of $\alpha_{n,k}=\check{\alpha}_{n,k}$ (note that $\hat{\alpha}_{n,k}=0$ in case $\theta=0$), condition \req{tangcone} and and minimality of $\rrr_l$.
Thus we have by $\omega_{n,k}\leq\overline{\omega}$ and Young's inequality that there exists $C>0$ such that
$$
\langle J_\pp^X(x_l-x_0)-J_\pp^X(x_m-x_0),x_l-\xsol\rangle_{X^*\times X}
\leq C\sum_{n=m}^{l-1} \left\{\left(\sum_{k=0}^{k_n-1} \ttt_{n,k}^\rr\right) + k_n \rrr_n^\rr\right)
$$
for which we can conclude convergence as $m,l\to\infty$ from \req{sumsumtnk} provided that
$$
\sum_{n=m}^\infty k_n \rrr_n^\rr \to 0 \mbox{ as } m\to\infty\,,
$$
which we guarantee by choosing, for an a priori fixed summable sequence $(a_n)_{n\in\N}$, e.g. $a_n=2^{-n}$
\begin{equation}\label{kn}
k_n := [a_n\rrr_n^{-\rr}]
\,,
\end{equation}
or, using \eqref{rnto0} just $k_n\equiv \bar{k}$ for some fixed integer $\bar{k}$, e.g., $\bar{k}=3$.
This is consistent with \eqref{alphanstar}, since in case $\delta=0$ we have $n_*=\infty$, so condition \eqref{alphanstar} never gets active in the noiseless case.


\subsection{Convergence with noisy data as $\delta\to0$}\label{subsec:conv_delta0}
In case $\nu,\theta>0$, convergence follows from the convergence rates results in Subsection \ref{subsec:rates}. Therefore it only remains to show convergence as $\delta\to0$ in case $\nu,\theta=0$.

In this section we explicitly emphasize dependence of the computed quantities on the noisy data and on the noise level by a superscript $\delta$.

Let $\|y^{\delta_j}-y\|\leq\delta_j$ with $\delta_j$ a zero sequence and $n_{*j}$ the corresponding stopping index. As usual \cite{HaNeSc95} we distinguish between the two cases that (i) $n_{*j}$ has a finite accumulation point and (ii) $n_{*j}$ tends to infinity.\\
\textbf{Case (i)}: there exists an $N\in\N$ and a subsequence $n_{j_i}$ such that for all $i\in\N$ we have $n_{j_i}=N$.
Provided
\begin{equation}\label{stability1}
n_*(\delta)=N\mbox{ for all }\delta \ \Rightarrow \
\mbox{The mapping }\delta\mapsto x_{N}^\delta\mbox{ is continuous at }\delta=0 \,,
\end{equation}
we can conclude that $x_N^{\delta_{j_i}}\to x_N^0$ as $i\to\infty$, and by taking the limit as $i\to\infty$ also in \req{discrprinc_upper}, $x_N^0$ is a solution to \req{Fxy}. Thus we may set $\xsol=x_N^0$ in \req{mon} (with $\theta=0$) to obtain
$$
\Dpx0{x_N^0}{z_{n_{*j_i},k_{n_{*j_i}}^*}^{\delta_{j_i}}}
=\Dpx0{x_N^0}{z_{N,k_{n_{*j_i}}^*}^{\delta_{j_i}}}
\leq \Dpx0{x_N^0}{x_N^{\delta_{j_i}}} \to 0 \mbox{ as }i\to\infty\,,
$$
where we have again used the continuous dependence \req{stability1} in the last step.\\
\textbf{Case (ii)}: let $n_{*j}\to\infty$ as $j\to\infty$, and let $\xsol$ be a solution to \req{Fxy}.
For arbitrary $\epsilon>0$, by convergence for $\delta=0$ (see the previous subsection) we can find $n$ such that
$\Dpx0{\xsol}{x_n^0}<\frac{\epsilon}{2}$ and, by Theorem 2.60 (d) in \cite{SKHK12}
there exists $j_0$ such that for all $j\geq j_0$ we have
$n_{*,j}\geq n+1$ and $|\Dpx0{\xsol}{x_n^{\delta_j}}-\Dpx0{\xsol}{x_n^0}|<\frac{\epsilon}{2}$, provided
\begin{equation}\label{stability2}
n\leq n_*(\delta)-1\mbox{ for all }\delta \ \Rightarrow \
\mbox{The mapping }\delta\mapsto x_n^\delta\mbox{ is continuous at }\delta=0 \,.
\end{equation}
Hence, by monotonicity of the errors we have
$$
\Dpx0{\xsol}{z_{n_{*j},k_{n_{*j}}^*}^{\delta_j}}
\leq \Dpx0{\xsol}{x_n^{\delta_j}}
\leq \Dpx0{\xsol}{x_n^0}+|\Dpx0{\xsol}{x_n^{\delta_j}}-\Dpx0{\xsol}{x_n^0}|<\epsilon\,.
$$

Indeed, \req{stability1}, \req{stability2} can be concluded from continuity of $F$, $F'$, the definition of the method \req{NIRLW}, as well as stable dependence of all parameters
$\omega_{n,k}$, $\alpha_{n,k}$, $k_n$ according to \req{omega_explicit}, \req{alphank0}, \req{alphank2}, \req{alphank02}, \req{kn} on the data $y^\delta$.

Altogether we have derived the following algorithm.
\begin{algorithm}\label{algoNIRLW} (Newton -- iteratively regularized Landweber method)
$$
\begin{array}{l}
\mbox{Choose $\Tau,\tilde{\tau},C_\alpha$ sufficiently large, $x_0$ sufficiently close to $\xsol$},\\
\mbox{$\alpha_{00}\leq 1$, $\vartheta>0$ sufficiently small, $\overline{\omega}>0$, $(a_n)_{n\in\N_0}$ such that $\sum_{n=0}^\infty a_n<\infty$}\\
\mbox{For $n=0,1,2\ldots$ until }\rrr_n\leq\Tau \delta\mbox{ do}\\
\hspace*{0.5cm}
u_{n,0}=0\\
\hspace*{0.5cm}
z_{n,0}=\xnd\\
\hspace*{0.5cm}
\alpha_{n,0} = \alpha_{n-1,k_{n-1}} \mbox{ if }n>0\\
\hspace*{0.5cm}
\mbox{For $k=0,1,2\ldots$ until } \left\{\begin{array}{ll}
k=k_n-1=a_n\rrr_n^{-\rr} &\mbox{ if } \rrr_n>\Tau \delta\\
\alpha_{n_*,k_{n_*}^*}\leq C_\alpha(\rrr_{n_*}+\delta)^{\frac{r}{1+\theta}}
&\mbox{ if }\rrr_n\leq\Tau \delta\end{array}\right\}
\mbox{ do}\\
\hspace*{1cm}
\begin{array}{l}
\omega_{n,k}=\vartheta\min\{
\ttt_{n,k}^{\frac{\rr}{\sd-1}}\tt_{n,k}^{-s}\,, \
\ttt_{n,k}^{\frac{\rr}{\pd-1}}\tt_{n,k}^{-p}\,, \ \overline{\omega}\}\\
u_{n,k+1}= u_{n,k}-\alpha_{n,k}J_\pp^X(z_{n,k}-x_0) \\
\qquad\qquad -\omega_{n,k} F'(\xnd)^*j_\rr^Y (F'(\xnd)(z_{n,k}-\xnd ) + F(\xnd)-y^\delta)\\
z_{n,k+1}=x_0+{J_\pp^X}^{-1} \Bigl(J_\pp^X(\xnd -x_0) + u_{n,k+1}\Bigr)\\
\alpha_{n,k+1}=\max\{\check{\alpha}_{n,k+1}\,, \, \hat{\alpha}_{n,k+1}\} \mbox{ with }
\check{\alpha}_{n,k+1}\,, \,\hat{\alpha}_{n,k+1} \mbox{ as in \req{alphank0}, \req{alphank2}}
\end{array}\\
\hspace*{0.5cm}
\xnpd=z_{n,k_n}
\,.
\end{array}
$$
Here we use the abbreviations according to \req{Anbn}, \req{ttt}, \req{tt}, \req{rrr}, \req{theta}.
\end{algorithm}

The analysis above yields the following convergence result.
\begin{theorem}\label{teo}
Assume that $X$ is smooth and $\ss$-convex with $\ss\geq\pp$,
that $x_0$ is sufficiently close to $\xsol$,
i.e., $x_0\in\mathcal{B}_{\rho}^D(\xsol)$,
that $F$ satisfies \req{tangcone} with \req{nonemptyint}, that
$F$ and  $F'$ are continuous and uniformly bounded in $\mathcal{B}_{\rho}^D(\xsol)$, and that \req{pnu}, \req{rp1} hold.
\\
Then, the iterates $z_{n,k}$ defined by Algorithm \ref{algoNIRLW} remain in $\mathcal{B}_{\rho}^D(\xsol)$ and converge to a solution $\xsol$ of \req{Fxy} subsequentially as $\delta\to0$ (i.e., there exists a convergent subsequence and the limit of every convergent subsequence is a solution).
In case of exact data $\delta=0$, we have subsequential convergence of $x_n$ to a solution of \req{Fxy} as $n\to\infty$.

If additionally a variational inequality \req{vinuNIRLW} with $\nu\in(0,1]$ and $\beta$ sufficiently small is satisfied, we obtain optimal convergence rates
\begin{equation}\label{ratesnuNIRLW}
\Dpx0{\xsol}{z_{n_*,k_{n_*}^*}}=O(\delta^{\frac{4\nu}{2\nu+1}})
\,,\quad \mbox{as} \;\; \delta \to 0\,.
\end{equation}

\end{theorem}

Note that we here deal with an a priori parameter choice: $\theta$ and therefore $\nu$ has to be known, otherwise $\nu$ must be set to a lower bound $\underline{\nu}$ for the true $\nu$ and since $\underline{\nu}\leq\nu$ implies validity of \req{vinuNIRLW} with $\nu$ replaced by $\underline{\nu}$, Theorem \ref{teo} still implies the (possibly suboptimal) rates $O(\delta^{\frac{4\underline{\nu}}{2\underline{\nu}+1}})$, or just convergence if we have set $\underline{\nu}=0$.

\section{Numerical Experiments}\label{sec:numres}
In this section we present some numerical experiments to test the method defined in section \ref{sec:conv_an}. We consider the identification of the space-dependent coefficient $c$ in the elliptic boundary value problem
\begin{equation}\label{boundaryproblem}
  \left\{
\begin{array}{ll}
-\Delta u+c u=f,& \text{ in } \Omega\\
u=0 & \text{ on }\partial \Omega\\
\end{array}
\right.
\end{equation}
from the measurement of $u$ in $\Omega$, where $f$ is a fixed function and where $\Omega$ is assumed to be a smooth, bounded domain in $\mathbb{R}^d$, $d$ $\in$ $\mathbb{N}$. Note that inhomogeneous Dirichlet boundary conditions can be easily incorporated into the right-hand side $f$ if necessary.\\
We consider three examples with $d=1$ and an example with $d=2$. In all cases, we take $X:=L^p(\Omega)$, $1<p<\infty$, $Y:=L^r(\Omega)$, $1<r<\infty$ and recall that the following facts hold true.
\begin{itemize}
  \item[(i)] For $1<p<\infty$, the duality mapping in $X$ is given by $J_p^X(c)=|c|^{p-1}\text{sgn}(c)$.
  \item[(ii)] If the domain of the forward operator is defined by
  \begin{equation}
    \mathcal{D}(F):=\{c \in X \text{ }|\text{ } \|c-\hat{c}\|_X \leq \gamma_p, \text{ for some } \hat{c} \in L^\infty(\Omega), \hat{c}\geq 0 \text{ a.e.}  \},
  \end{equation}
  the condition \req{nonemptyint} is satisfied with $c_0=0$.
  \item[(iii)] There follows from Lemma $2$ in \cite{bkSchoepferSchuster09} that the operator $F: \mathcal{D}(F) \subseteq X \rightarrow Y$, $F(c)=A(c)^{-1}f$ is well defined. Here $A(c): W^{2,p}(\Omega) \cap H_0^{1}(\Omega) \rightarrow L^p(\Omega)$ is given by $A(c)u:=-\Delta u+cu$. Moreover, $F'(c)$ and the adjoint of $F'(c)$
   $$F'(c)h=-A(c)^{-1}(h F(c)), \text{ }\text{ }\text{ }F'(c)^*w=-u(c) A(c)^{-1}w,$$
   are well defined and bounded.
\end{itemize}
In all the numerical simulations, we take $\theta=0$ and stop the outer iteration by means of the discrepancy principle \eqref{discrprinc_upper}.
Concerning the stopping index of the inner iteration, we slightly modify Algorithm 1, requiring also that if $\|F(z_{n,k})-y^\delta\| \leq \tau \delta$ then the iteration has to be stopped. More precisely,
\begin{equation}\label{knnumexp}
  k_n=\min \{k \in \mathbb{Z}, k \geq 0,\text{ } \text{ } | \text{ }  \|F(z_{n,k})-y^\delta\|\leq \tau \delta \text{ } \vee \text{ } k \geq a_n \rrr_n^{-r}\}
\end{equation}
and the regularized solution is $c_{n_*}^\delta=z_{n_{n_*-1},k_{n_*-1}}$.

\subsection{1-dimensional examples}
We consider the same numerical simulations as in \cite{KaltTomb13}, 
taking $\Omega=(0,1)$ and inhomogeneous boundary conditions $u(0)=g_0$, $u(1)=g_1$. We solve all differential equations approximately by a finite difference method by dividing the interval $[0,1]$ into $N+1$ subintervals with equal length $1/(N + 1)$, in all examples below $N = 400$. The $L^p$ and $L^r$ norms are calculated approximately by means of a quadrature method.

\begin{example}\label{doublepeak}
In the first simulation we assume that the solution is sparse:
\begin{equation}\label{cdagdef}
c^\dag(t)=
  \left\{
\begin{array}{ll}
0.5,  & 0.3 \leq t \leq 0.4,\\
1.0,  & 0.6 \leq t \leq 0.7,\\
0.0,  &elsewhere.\\
\end{array}
\right.
\end{equation}
The test problem is constructed by taking $u(t)=u(c^\dag)(t)=1+5t$, $f(t)=u(t)c^\dag(t)$, $g_0=1$ and $g_1=6$. We perturb the exact data $u$ with gaussian white noise: the corresponding perturbed data $u^\delta$ satisfies $\|u^\delta-u\|_{L^r}=\delta$, with $\delta=0.1 \times 10^{-3}$.\\
We apply Algorithm 1, with the inner stopping index satisfying \req{knnumexp}, with $\tau=1.02$, $\tilde{\tau}=0.1$, $a_n=(50+n)^{-2}$. The upper bound $\overline{c}_\omega$ is fixed equal to $0.1$ and $\vartheta$ is chosen as $2^{-j^\sharp}$, where $j^\sharp$ is the first index that satisfies
\begin{equation}
2^{s^*-1} C_{p^*,s^*}(p\rho^2)^{1-s^*/p^*}\vartheta^{s^*-1}+2^{p^*-1} C_{p^*,s^*} \vartheta^{p^*-1} \leq \overline{c}_\omega.
\end{equation}
\begin{figure}
\centering
\subfigure[$p=2$, $r=2$.]
{\includegraphics[width=6cm]{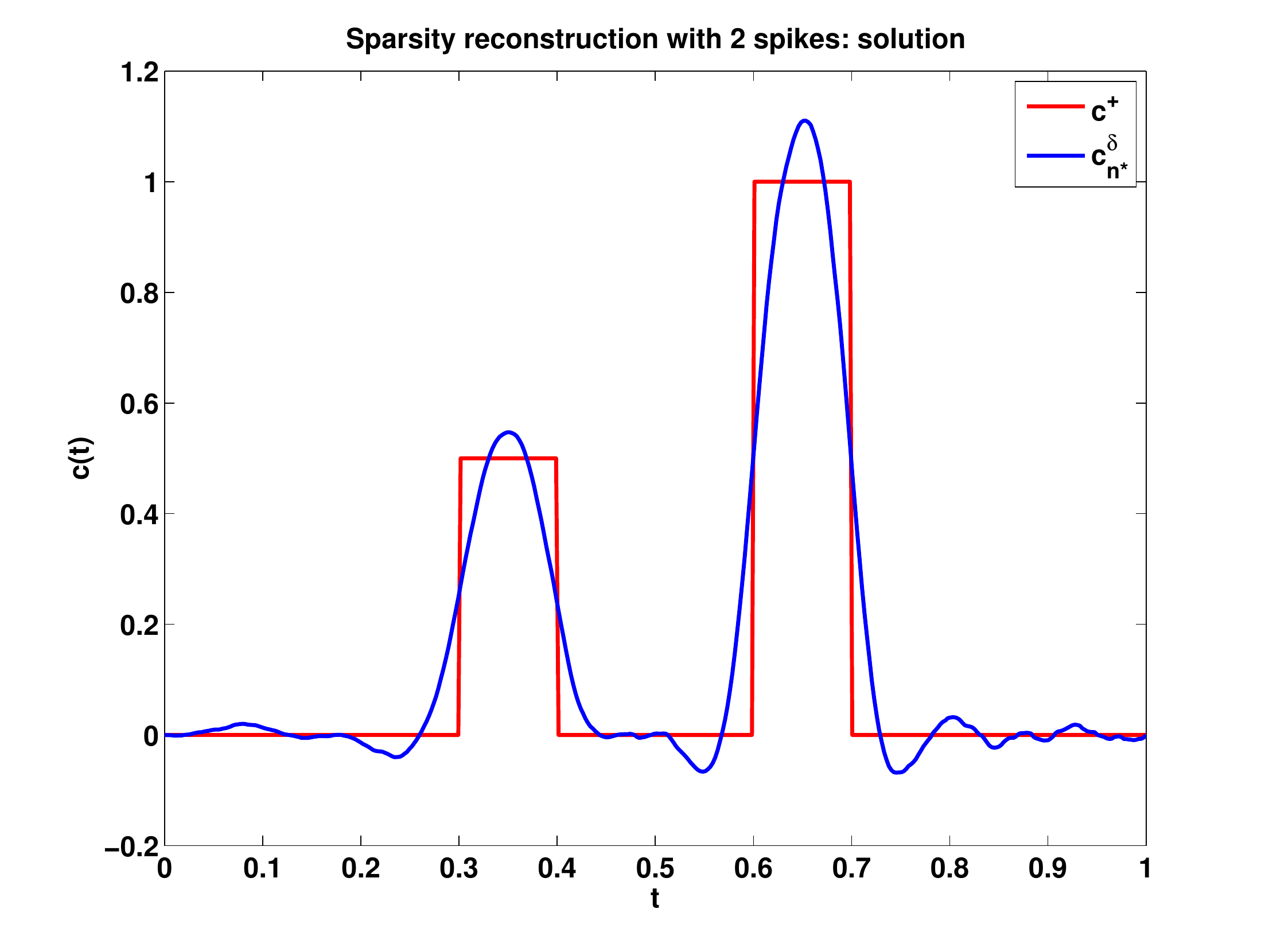}}
\subfigure[$p=1.1$, $r=2$.]
{\includegraphics[width=6cm]{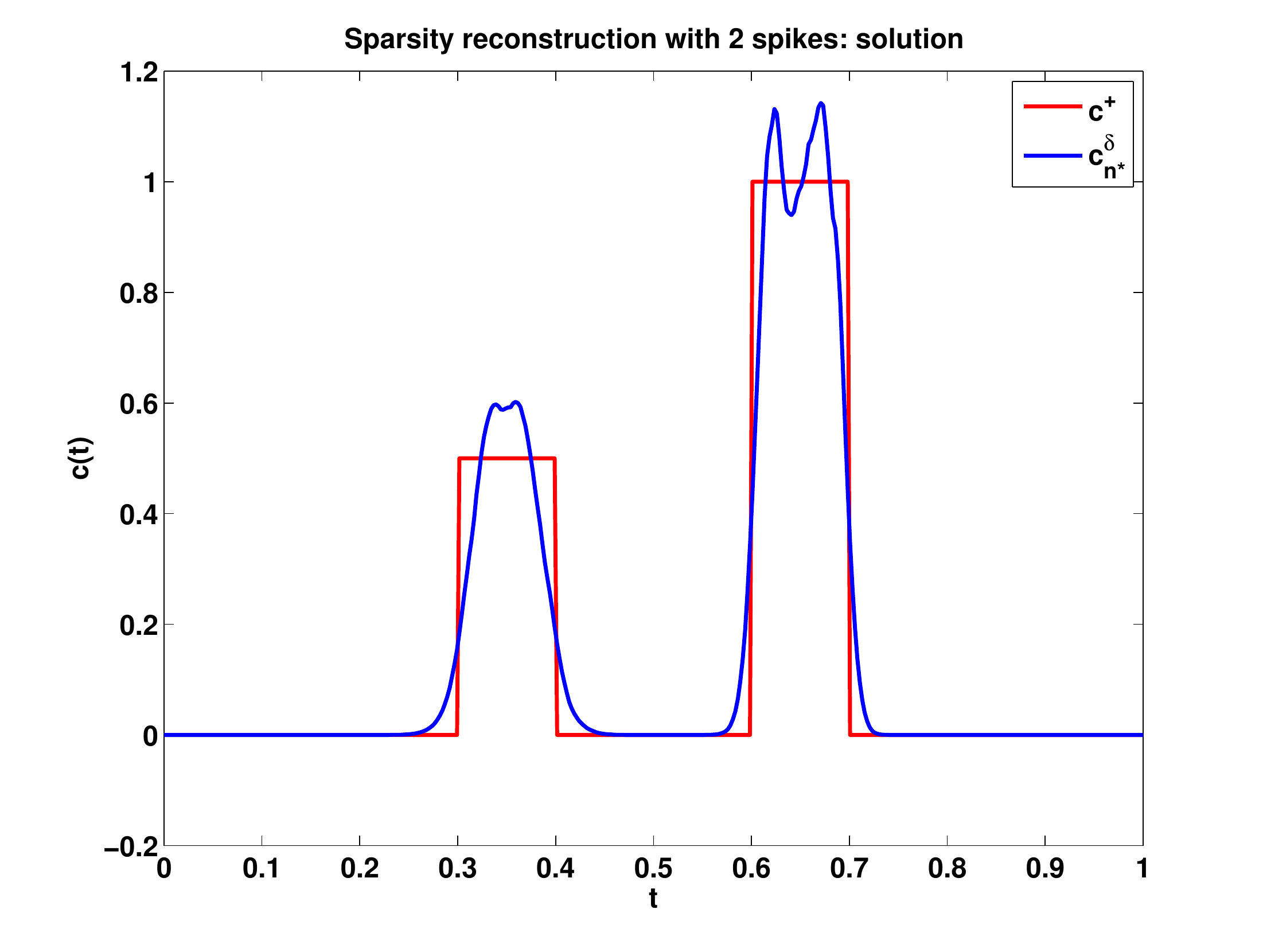}}
\caption{Reconstructed Solutions for example \ref{doublepeak}}
\label{figure1}
\end{figure}
In figure \ref{figure1} we show the results obtained by our method with $p=2$ and $p=1.1$ respectively. The reconstructed solutions are very similar to those obtained in \cite{KaltTomb13} and \cite{JinNLW}. Concerning the total number of inner iterations
$$N_p=\sum_{n=0}^{n_*-1} k_n,$$
similarly to \cite{KaltTomb13}, it is larger in the case $p=2$ ($N_2=3992$) than in the case $p=1.1$ ($N_{1.1}=3063$). In both cases, the value of $N_p$ is lower than the corresponding value found in \cite{KaltTomb13}. We underline that it is possible to make different choices for $\tilde{\tau}$, $\overline{c}_\omega$, $a_n$ to look for further improvements of the speed of the method. The partial freedom in the choices of these parameters makes Algorithm 1 more flexible than the method described in \cite{KaltTomb13}. In our numerical simulations, we tested different choices which gave similar but slightly worse results than those stated here.
\end{example}

\begin{example}\label{triplepeak}
We modify the exact solution of the previous example into:
\begin{equation}\label{cdagdef2}
c^\dag(t)=
  \left\{
\begin{array}{ll}
0.25, & 0.1  \leq t \leq 0.15,\\
0.5,  & 0.3  \leq t \leq 0.4,\\
1.0,  & 0.6  \leq t \leq 0.7,\\
0.0   &elsewhere.\\
\end{array}
\right.
\end{equation}
and choose again $\delta=0.1 \times 10^{-3}$. In this case, we take $\tau=1.02$, $\tilde{\tau}=0.01$, $a_n=(100+n)^{-2}$, $\overline{c}_\omega=0.1$ and $\vartheta$ as in the previous example.
\begin{figure}
\centering
\subfigure[$p=2$, $r=2$.]
{\includegraphics[width=6cm]{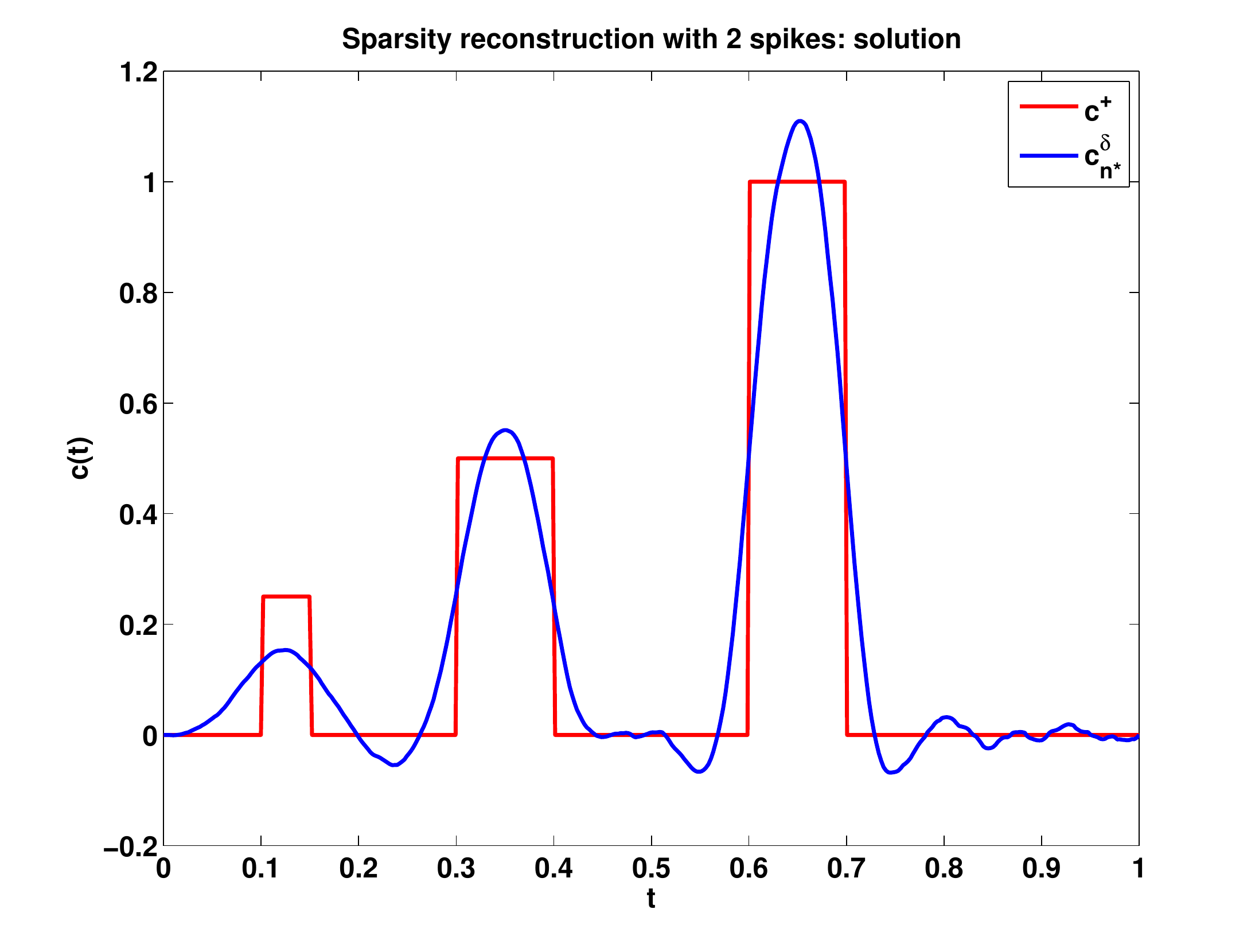}}
\subfigure[$p=1.1$, $r=2$.]
{\includegraphics[width=6cm]{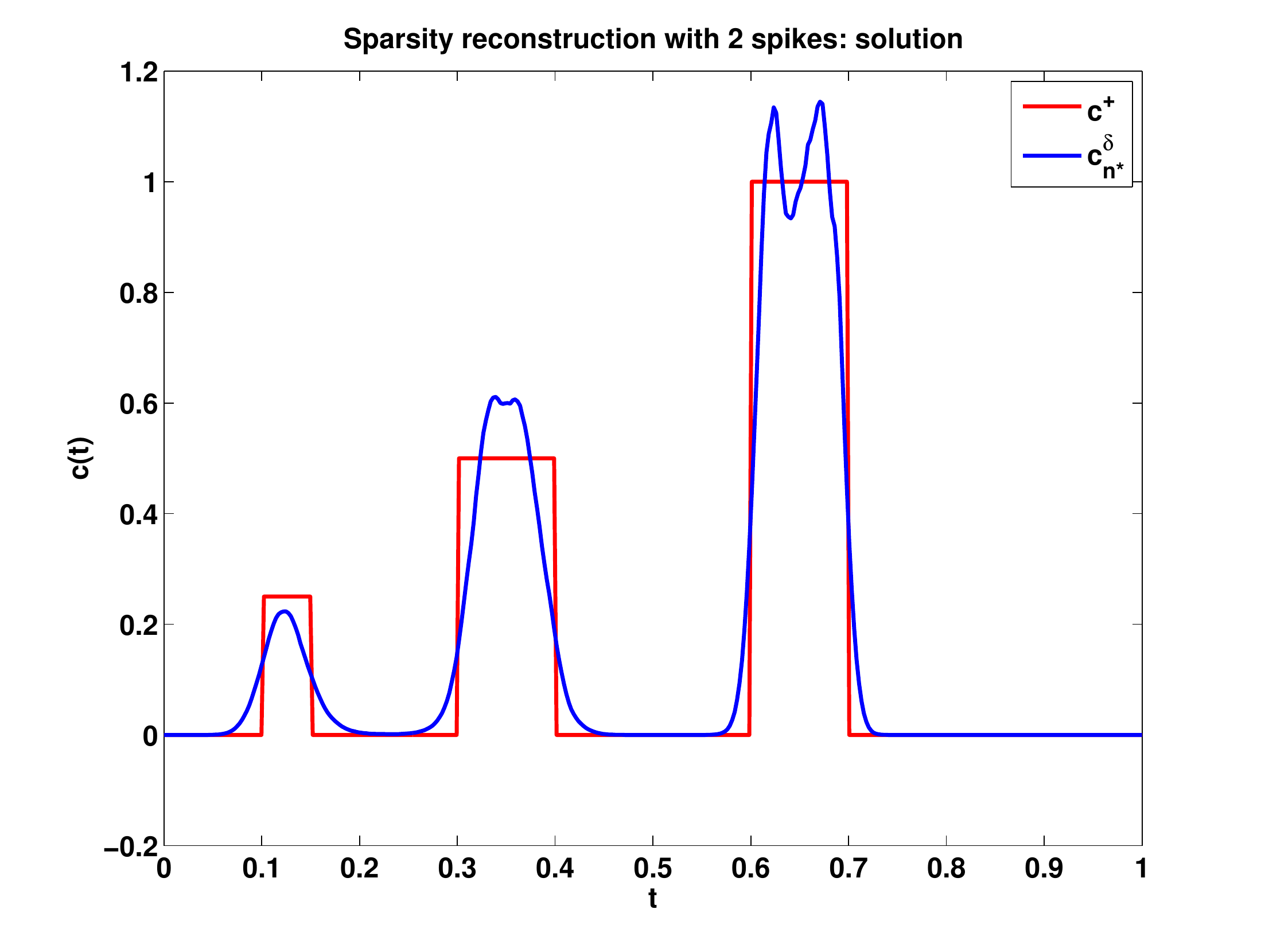}}
\caption{Reconstructed Solutions for example \ref{triplepeak}}
\label{figure2}
\end{figure}
In figure \ref{figure2} we show the results obtained by our method with $p=2$ and $p=1.1$ respectively. As usual, the reconstruction of the sparsity is much better for $p=1.1$. Concerning the total number of inner iteration $N_p$, in this case we obtain $N_2=4141$ (with a corresponding error of $0.1110$) and $N_{1.1}=3110$ (with a corresponding error of $0.0482$). We observe that although the corresponding error is slightly larger than that obtained in \cite{KaltTomb13}, the value of $N_2$ is slightly more than a fifth than the value obtained in \cite{KaltTomb13}, with a gain in the speed of the $421.8\%$. Moreover, in the case $p=1.1$ Algorithm 1 performs $1415$ iterations less than in \cite{KaltTomb13}, with a gain in the speed of the $45.5\%$, and obtains even a lower error.
\end{example}

\begin{example}\label{Outliers}
We consider an example with noisy data where a few data points called outliers are remarkably different from other data points. This situation may arise from procedural measurement errors.\\
We suppose $c^\dag$ to be a smooth solution
\begin{equation}
  c^\dag(t)=2-t+4\sin(2\pi t)
\end{equation}
and take $u(c^\dag)(t)=1-2t$, $f(t)=(1-2t)(2-t+4\sin(2\pi t))$, $u(0)=g_0=1$ and $u(1)=g_1=-1$ as exact data of the problem. We start the iteration from the initial guess $c_0(t)=2-t$, and take $\tilde{\tau}=5\times 10^{-3}$, $a_n=(1+n)^{-1.1}$, $\overline{c}_\omega=5\times 10^{-3}$ and $\vartheta$ as in the previous example. Using the same Matlab seed for generating random data, we consider the same perturbed data as in \cite{KaltTomb13}, Example 3, case $B$ (cf. Figure \ref{figure3} here and Figure 3, picture (d) in that paper).\\
We run both Algorithm $1$ with $\tilde{\tau}=1.0015$ and the algorithm that generated the results presented in \cite{KaltTomb13} for this example with $p=2$ and $r=1.1$\footnote{These values of $r$ and $p$ do not satisfy the condition \req{pnu} of Theorem \ref{teo}. However, this condition is needed only to have a lower bound for $\omega_{n,k}$ and such a bound is verified experimentally in this example for these values of $r$ and $p$.}.
\begin{figure}
\centering
\subfigure[Perturbed data with outliers]
{\includegraphics[width=6cm]{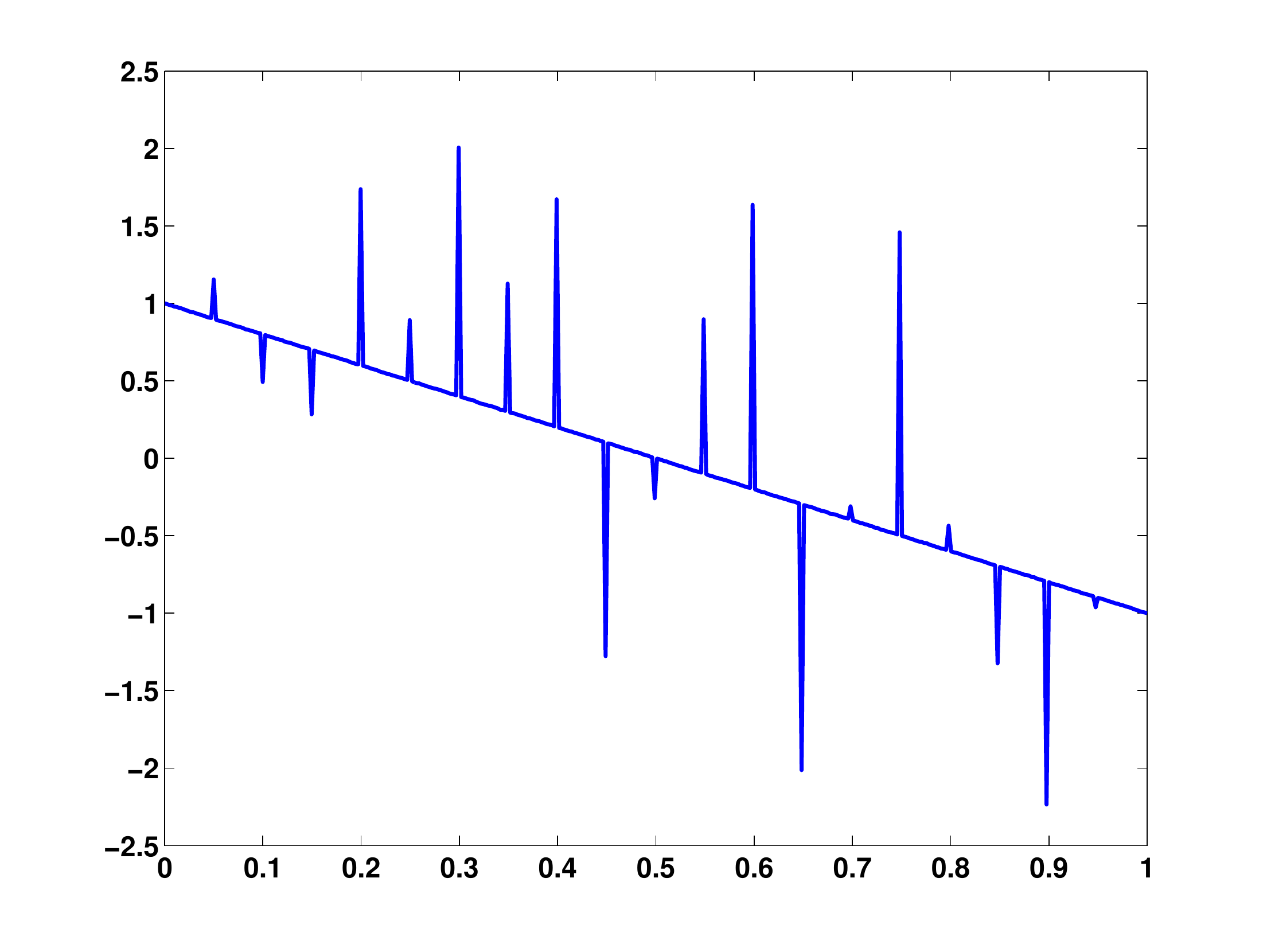}}
\subfigure[Alg. 1 Solution with $\tau=1.0015$]
{\includegraphics[width=6cm]{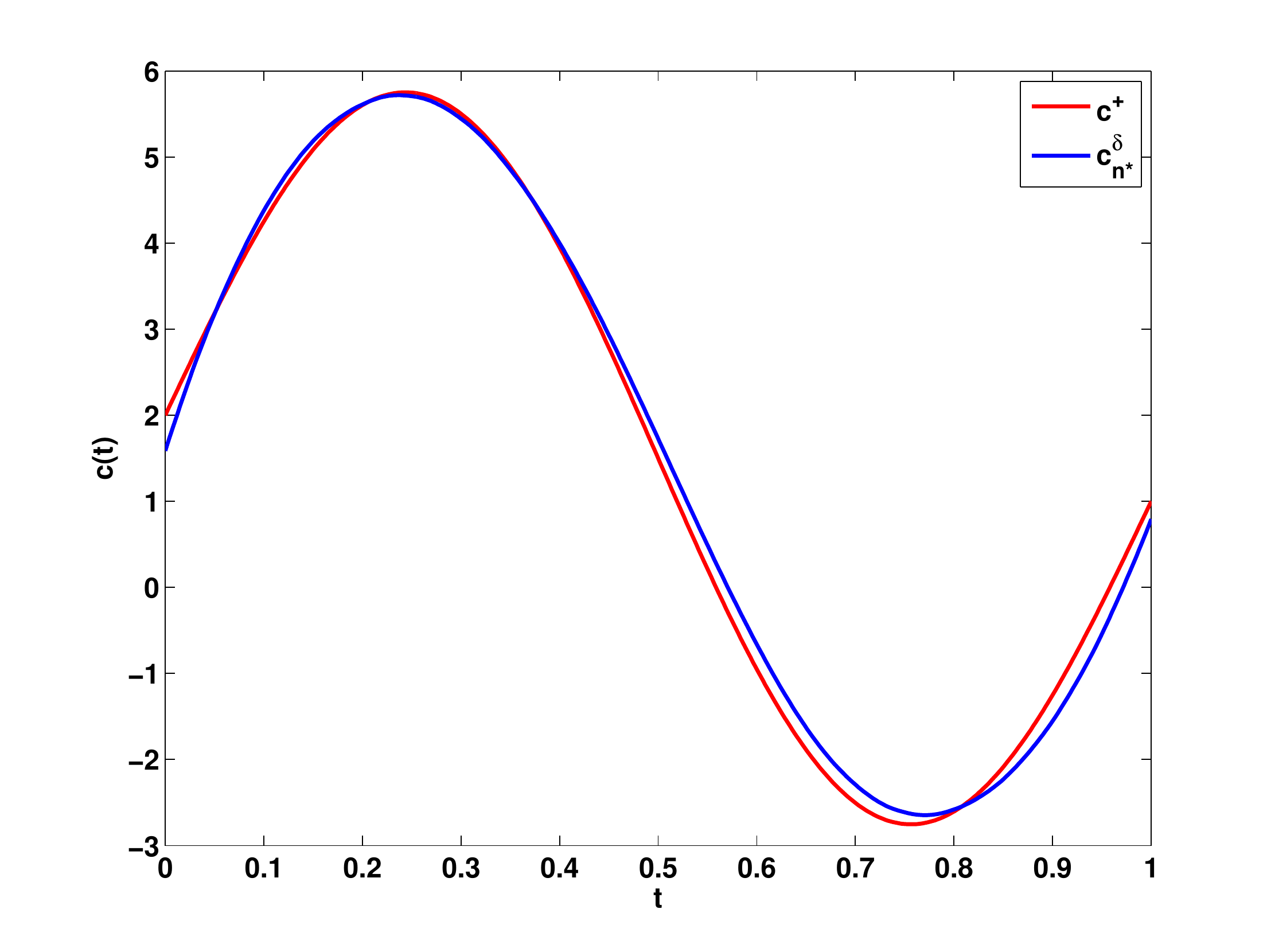}}
\subfigure[Alg. 1 Solution with $\tau=1+10^{-5}$]
{\includegraphics[width=6cm]{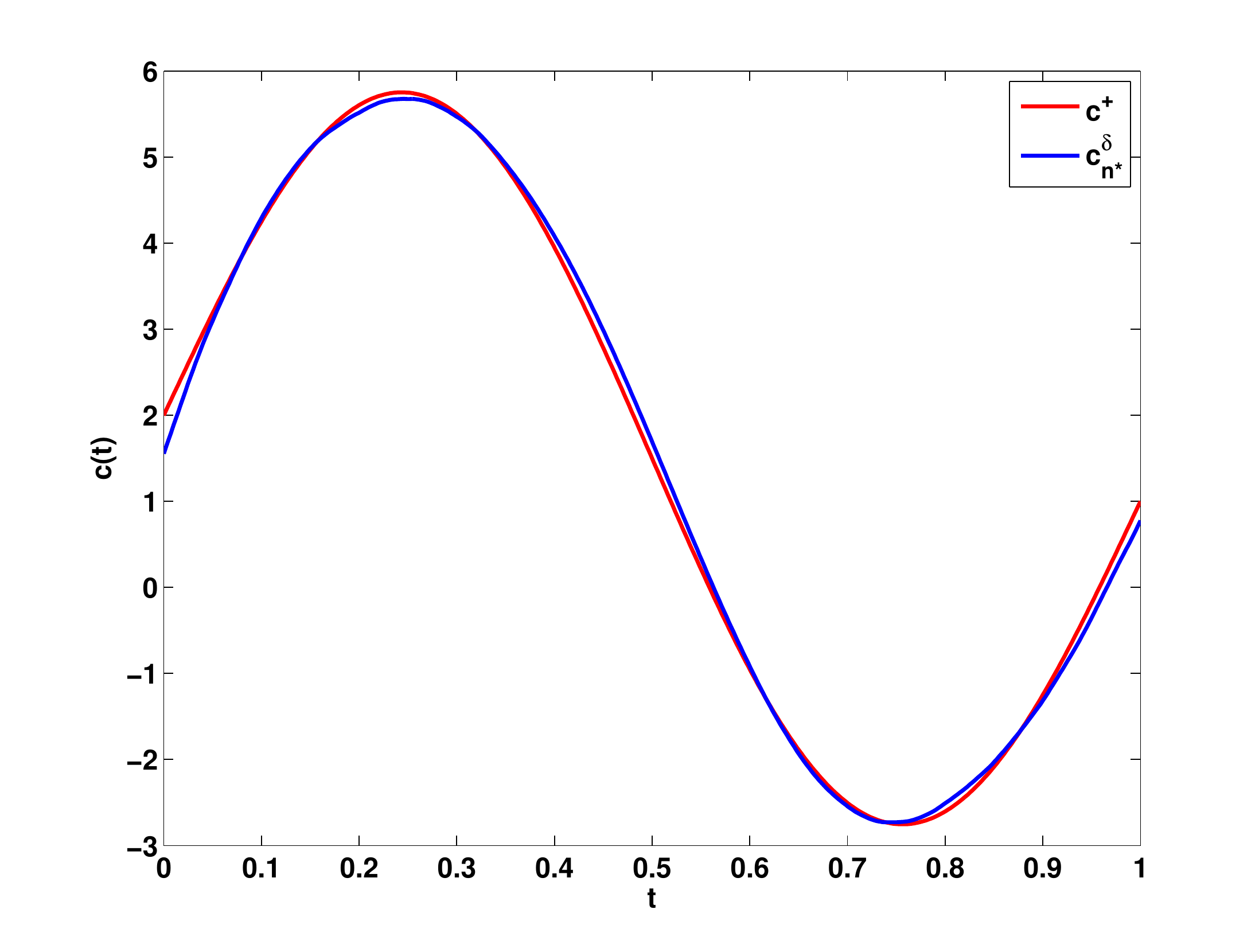}}
\subfigure[Solution of the method in \cite{KaltTomb13}]
{\includegraphics[width=6cm]{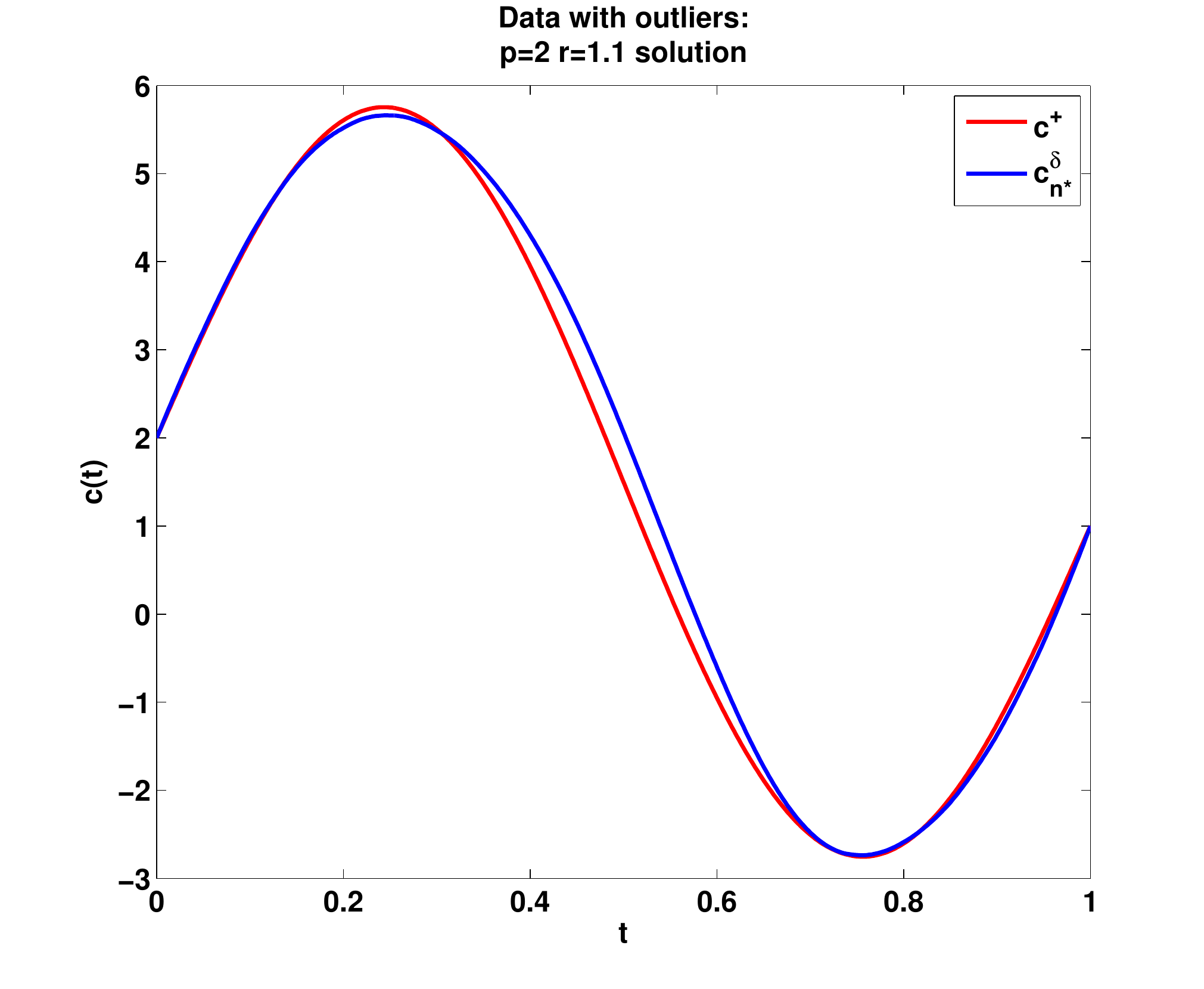}}
\caption{Numerical results for example \ref{Outliers}: (a) are the perturbed data; (b) is the solution obtained by the method in \cite{KaltTomb13}; (c) and (d) are the solutions of Algorithm 1 with different values of $\tau$.}
\label{figure3}
\end{figure}
Pictures (b) and (d) from figure \ref{figure3} show the corresponding results. The solution obtained by Algorithm 1 is slightly more precise, with an error equal to $1.8852 \times 10^{-1}$ for Algorithm 1 and equal to $2.9388 \times 10^{-1}$ for the method in \cite{KaltTomb13}. The most interesting fact is that Algorithm $1$ computes only $N_2=249$ total inner iterations to obtain this solution, whereas for the method in \cite{KaltTomb13} $N_2=3285$ and the reconstruction is poorer. Moreover, due to the flexibility of our method, we can simply change the value of $\tau$ into $1+10^{-5}$ to get a more precise solution (see picture (c) in figure \ref{figure3}). A visual inspection gives an idea of the improvement: the error of this solution, obtained with $278$ total inner iterations, is equal to $1.1607 \times 10^{-1}$.\\
We summarize the numerical results of the 1-dimensional examples in Table $\ref{ResultsExample2}$.
\begin{table}[h]
\centering
\begin{tabular}{|c|c|c|}\hline
\multicolumn{3}{|c|}{\textbf{1-dimensional numerical simulations}}\\ \hline
\multicolumn{3}{|c|}{\textbf{Example 1}}\\ \hline
Method & Total n. of inner iterations & Error: $\|\cdot-c^\dag\|$\\ \hline
Alg. 1 $p=1.1$            & 3063       & 0.0413\\
Alg. 1 $p=2.0$            & 3992       & 0.1059\\ \hline
\multicolumn{3}{|c|}{\textbf{Example 2}}\\ \hline
Method & Total n. of inner iterations & Error: $\|\cdot-c^\dag\|$\\ \hline
Alg. 1 $p=1.1$            & 3110       & 0.0482\\
Alg. 1 $p=2.0$            & 4141       & 0.1110\\  \hline
\multicolumn{3}{|c|}{\textbf{Example 3}}\\ \hline
Method & Total n. of inner iterations & Error: $\|\cdot-c^\dag\|$\\ \hline
Alg. 1 $\tau=1.0015$               & 249        & 0.1885\\
Alg. 1 $\tau=1+10^{-5}$            & 278        & 0.1161\\
Method from \cite{KaltTomb13}       & 3285       & 0.2939\\ \hline
\end{tabular}
\caption{Numerical results for the 1-dimensional simulations.}
\label{ResultsExample2}
\end{table}
\end{example}
\subsection{A 2-dimensional example}
We show the performance of Algorithm 1 in the following 2-dimensional example.\\
Let $\Omega=(0,1)^2$ and assume the exact solution to be
\begin{equation}\label{cdagdef4}
c^\dag(x,y)=40 \chi_{[0.19,0.24]^2}(x,y), \text{ }\text{ }\text{ }x,y \in \Omega,
\end{equation}
where the function $\chi$ is the characteristic function on a subset of $\mathbb{R}^2$.
We take as exact data $u(x,y)=1+x+y$: as a consequence, the fixed right hand side of the problem is $f(x,y)=c^\dag(x,y)u(x,y)$ and the data at the boundary are given by $g(x,y)=u(x,y)$ for every $(x,y)$ $\in$ $\partial$ $\Omega$. We discretize the interval $[0,1]_x$ into $N+1$ subintervals and the interval $[0,1]_y$ into $M+1$ subintervals and compute the solutions of the forward operator $F(c)=u(c)$ by a finite difference method. The $L^p$ and $L^r$ norms in $\Omega$ are calculated by a simple quadrature method.\\
We fix $p=1.1$, $N=M=30$, $\tau=1+10^{-5}$, $\tilde{\tau}=10^{-4}$, $\overline{c}_\omega=0.1$ and $\vartheta$ as in the 1-dimensional examples. We run Algorithm 1 starting from $c_0=0$ with the data $u^\delta$ perturbed by white gaussian noise with two different noise levels. In the first case, we choose a small value for $\delta=\|u-u^\delta\|=10^{-3}$ with $r=2$. In the second case we choose $\delta=10^{-2}$ with $r=2$ and with $r=10$.\\
\begin{figure}
\centering
\subfigure[Exact solution]
{\includegraphics[width=6cm]{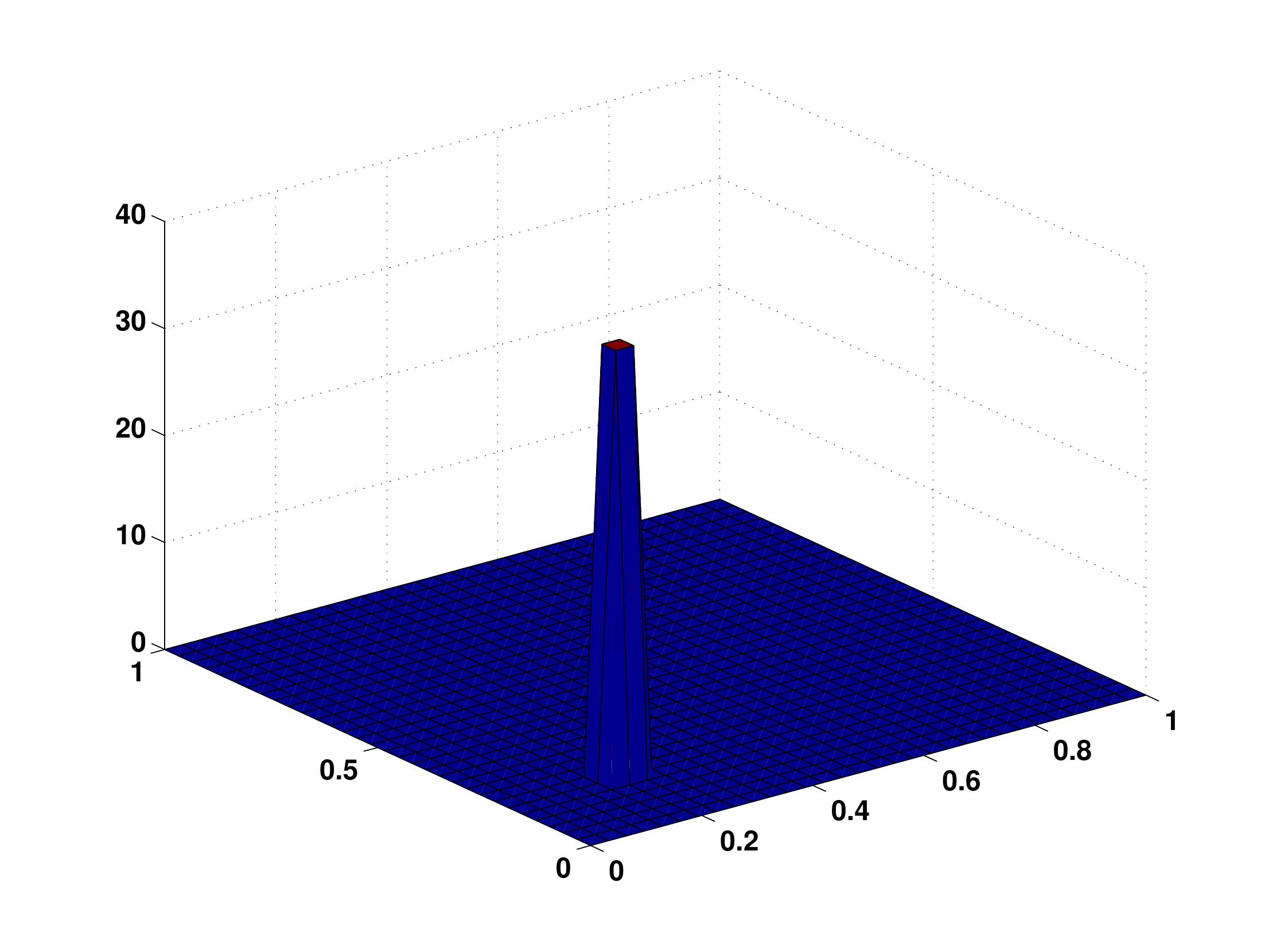}}
\subfigure[]
{\includegraphics[width=6cm]{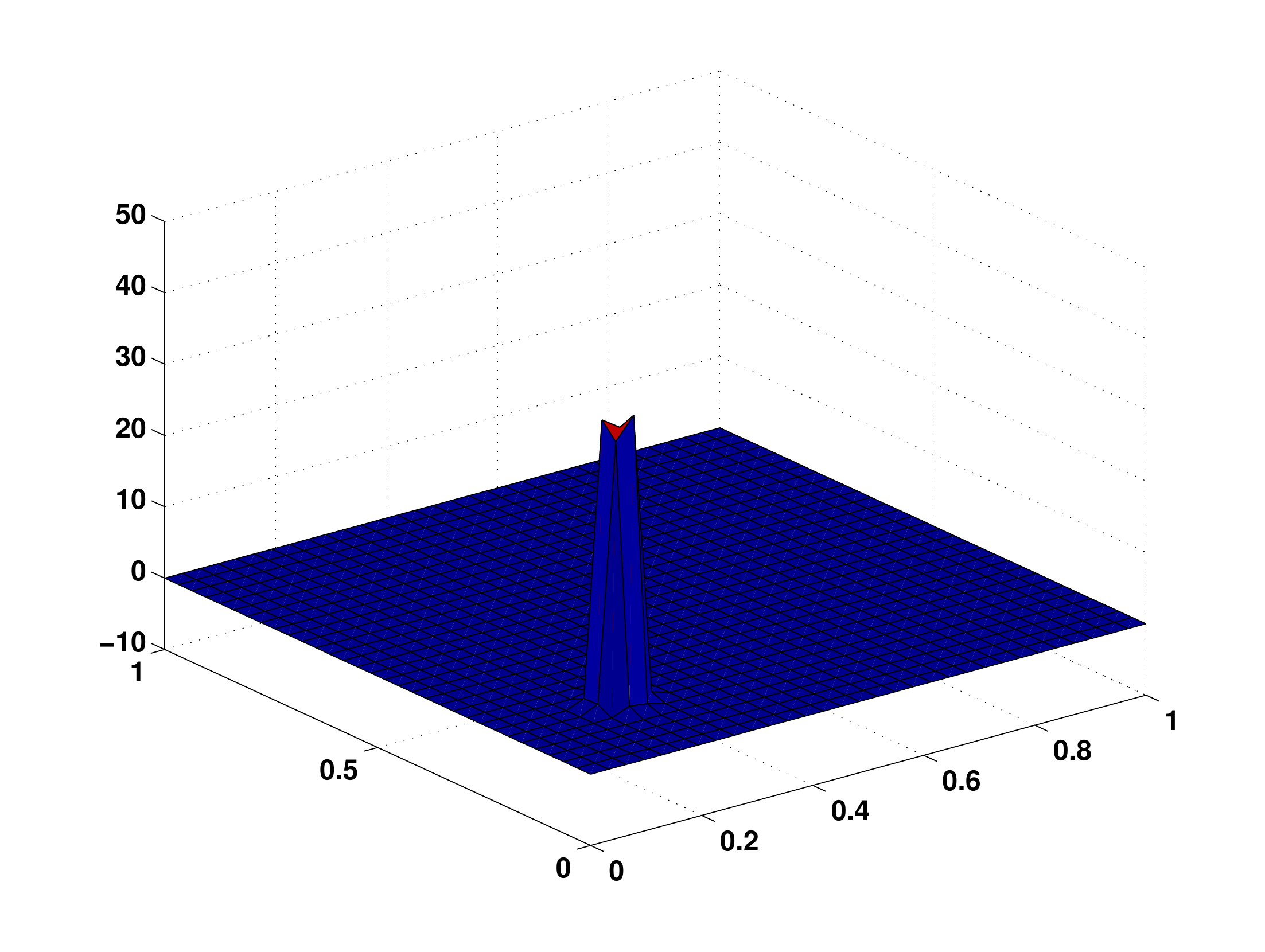}}
\subfigure[]
{\includegraphics[width=6cm]{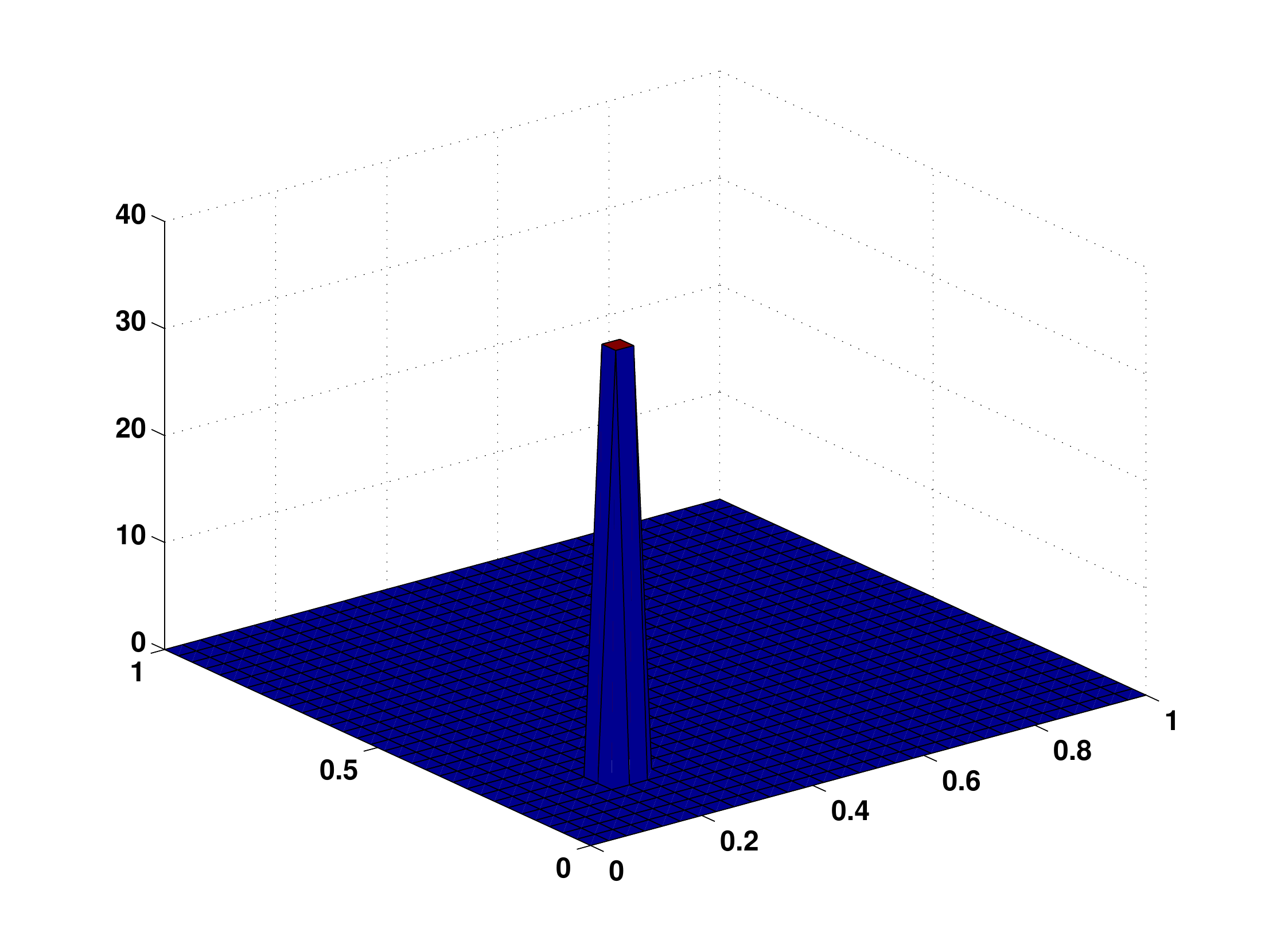}}
\subfigure[]
{\includegraphics[width=6cm]{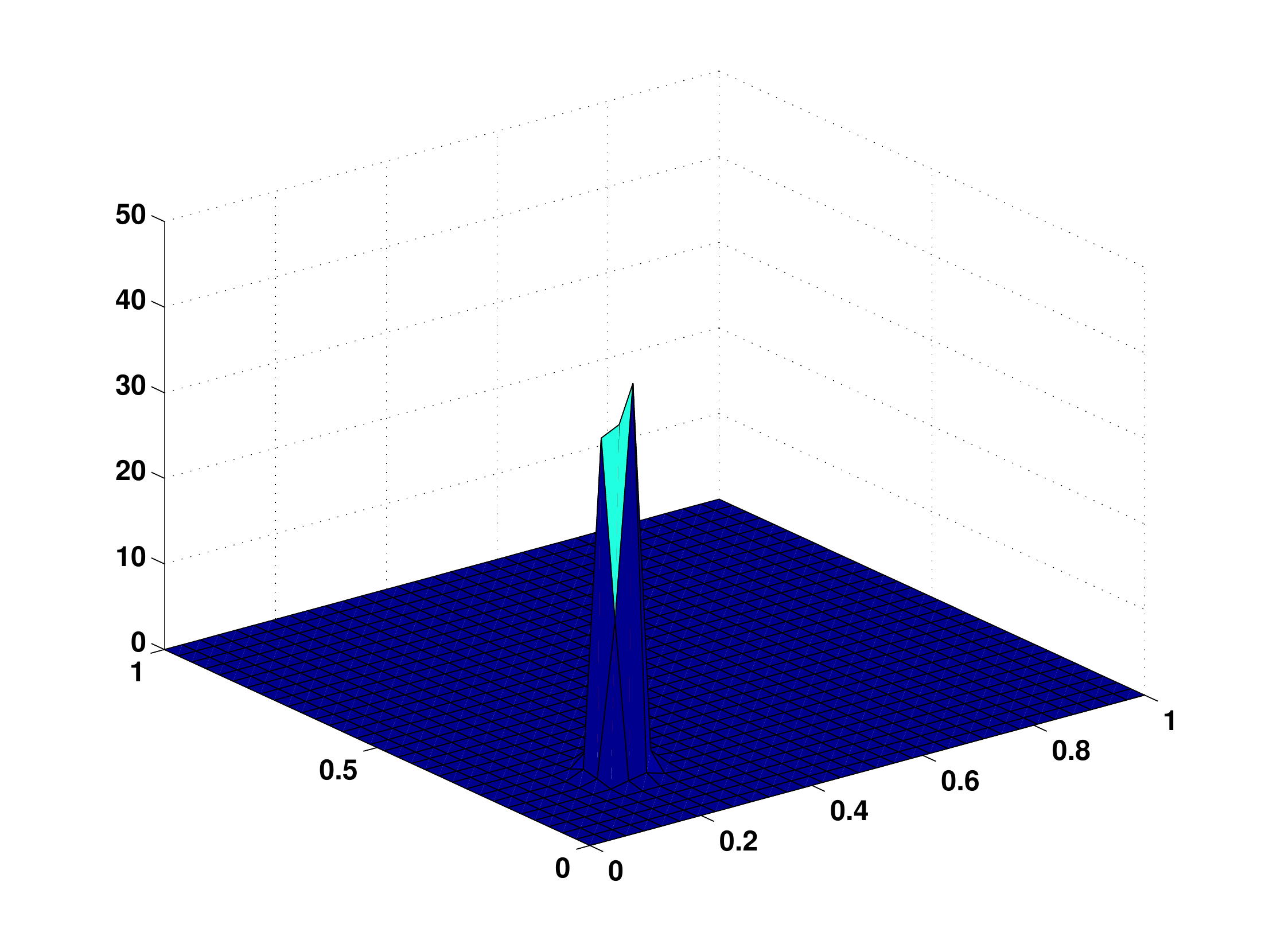}}
\caption{Numerical results for the 2-dimensional example: (a) is the exact solution; (b) is the solution corresponding to $\delta=10^{-3}$ and $r=2$; (c) is the solution corresponding to $\delta=10^{-2}$ and $r=2$; (d) is the solution corresponding to $\delta=10^{-2}$ and $r=10$.}
\label{figure4}
\end{figure}
In figure \ref{figure4} we plot the exact solution and the reconstructions obtained using Algorithm 1 in all these cases. The pictures show that the method provides a good reconstruction of the sparsity in this example. In particular, we underline that in the case $\delta=10^{-2}$ the choice of a large $r$ improves the result, obtaining a better reconstruction with very few iterations (the total number of iterations is equal to $9$ in this case).

\section{Conclusions}
In this paper we have devised an alternative parameter choice strategy for the iteratively regularized Newton- Landweber iteration proposed in \cite{KaltTomb13}. This strategy is based on alternative error estimates and allows for a unified treatment of the unconditional convergence case and convergence with rates.
In our future research we will try to extend the analysis to faster inner iterations than Landweber such as steepest descent or conjugate gradient methods.
\section{Acknowledgment}
Support by the German Science Foundation DFG under grant KA 1778/5-1 is gratefully acknowledged.

\end{document}